\documentclass[11pt,reqno]{amsart}
\usepackage[T1]{fontenc}
\usepackage[utf8]{inputenc}
\usepackage{lmodern}
\usepackage[a4paper,textwidth=158mm,textheight=240mm,centering]{geometry}
\usepackage{amssymb,mathtools,mathrsfs,booktabs}
\usepackage{microtype}
\usepackage{enumitem}
\usepackage[hidelinks,unicode]{hyperref}
\usepackage{bookmark}
\numberwithin{equation}{section}
\allowdisplaybreaks[2]
\setlist[enumerate]{label=\textup{(\roman*)},leftmargin=2.15em,itemsep=0.25em,topsep=0.4em}
\theoremstyle{plain}
\newtheorem{theorem}{Theorem}[section]
\newtheorem{proposition}[theorem]{Proposition}
\newtheorem{lemma}[theorem]{Lemma}
\newtheorem{corollary}[theorem]{Corollary}
\theoremstyle{definition}

\newtheorem{example}[theorem]{Example}
\newtheorem{conjecture}[theorem]{Conjecture}
\theoremstyle{remark}
\newtheorem{remark}[theorem]{Remark}
\newcommand{\R}{\mathbb R}

\newcommand{\A}{\mathscr A}
\newcommand{\B}{\mathcal B}
\newcommand{\cL}{\mathcal L}
\newcommand{\cP}{\mathcal P}
\newcommand{\cT}{\mathcal T}
\DeclareMathOperator{\Sym}{Sym}
\DeclareMathOperator{\rank}{rank}
\DeclareMathOperator{\tr}{tr}
\DeclareMathOperator{\im}{im}
\DeclareMathOperator{\Span}{span}
\DeclareMathOperator{\aff}{aff}
\DeclareMathOperator{\dist}{dist}
\DeclareMathOperator{\proj}{proj}
\DeclareMathOperator{\Graff}{Graff}

\DeclareMathOperator{\adj}{adj}
\DeclareMathOperator{\diag}{diag}
\DeclareMathOperator{\Syl}{Syl}
\newcommand{\T}{^{\mathsf T}}
\newcommand{\norm}[1]{\left\lVert#1\right\rVert}
\newcommand{\ip}[2]{\left\langle#1,#2\right\rangle}
\title[Equidistant affine subspaces]{Equidistant Affine Subspaces:\\Explicit Bounds, Structural Reductions,\\Transverse Products, and Fixed-Codimension Stability}
\author{Jia Li}
\address{School of Mathematical Sciences,  Key Laboratory of MEA (Ministry of Education) \& Shanghai Key Laboratory of PMMP,  East China Normal University, Shanghai, China}
\email{jialimath001@pku.org.cn}
\date{}
\hypersetup{pdftitle={Equidistant Affine Subspaces: Explicit Bounds, Structural Reductions, Transverse Products, and Fixed-Codimension Stability},pdfauthor={Jia Li}}
\begin{document}
\begin{abstract}
Let $\mathcal{N}(k,n)$ be the largest size of a family of affine $k$-planes
in $\R^n$ at pairwise separation distance one. We develop a unified
treatment of upper bounds, dimension reduction, product constructions,
and fixed-codimension stability. Orthogonal-projection coordinates give a
rank-stratified distance identity; a Veronese lifting converts its
polynomial consequences into simultaneous total isotropy. Connected Nash
decomposition and an analytic diagonal argument yield an explicit bound
\[
\mathcal{N}(k,n)\leqslant 2^{\,2^{\,3(k+1)(n-k)}}.
\]
For configurations whose direction sums have dimension at most $k+1$,
we prove the exact reduction
$\mathcal{N}_{\leqslant k+1}(k,n)=\mathcal{N}(1,n-k+1)$.
Every codimension-two configuration splits off a common linear factor,
so $\mathcal{N}(n-2,n)=\mathcal{N}(1,3)$. A transverse product construction
gives multiplicative inequalities. We further prove that for every fixed
codimension $c$, the sequence $\mathcal{N}(k,k+c)$ is eventually constant.
In codimension three, a sharper normal-span argument gives
\[
\mathcal{N}(k,k+3)=\mathcal{N}(7,10)\qquad(k\geqslant7),
\]
where the numerical threshold uses the computer-assisted theorem
$\mathcal{N}(1,3)=7$. The structural and stability results are otherwise
independent of that numerical input.
\end{abstract}
\maketitle

\section{Introduction}\label{sec:introduction}
For affine subspaces $L,L'\subset\R^n$, their \emph{separation distance} is
\begin{equation}\label{eq:distance-definition}
 d(L,L')=\inf\left\{\norm{x-y}\big|x\in L,\ y\in L'\right\}.
\end{equation}
For $k=0$, this is the usual Euclidean distance between points. For $k>0$,
this is not a metric on the affine Grassmannian and generally differs from
the Hausdorff distance. For integers $0\leqslant k<n$, define
\begin{equation}\label{eq:N-definition}
 \mathcal{N}(k,n):=\sup\left\{|\cL|\bigg|
 \begin{array}{l}
 \cL\subset\Graff(k,n)\text{ is finite},\\
 d(L,L')=1\text{ whenever }L\ne L'\text{ in }\cL
 \end{array}\right\}.
\end{equation}
Initially the value $+\infty$ is allowed. A uniform bound on all finite
families also excludes infinite families. Throughout, an equidistant family
means a family with a common strictly positive separation distance,
normalized to one. Rescaling shows that the same extremal number applies
to any prescribed common positive distance.
Throughout the paper, $k$ is the dimension of each affine subspace and
$n$ is the ambient dimension.

For points, the regular simplex gives $\mathcal{N}(0,n)=n+1$; for hyperplanes,
$\mathcal{N}(n-1,n)=2$. The case $\mathcal{N}(1,3)$ is the axis formulation of Littlewood's
problem on mutually touching congruent infinite cylinders
\cite{Littlewood}. Boz\'oki, Lee, and R\'onyai established the existence of
seven such cylinders \cite{BLR}. H\"ofer proved $\mathcal{N}(1,3)\leqslant9$ by a
computer-free linear-algebraic argument \cite{Hoefer}, following the
Ramsey-theoretic bound ten of Dillon, Koizumi, and Luo \cite{DKL}.
A recent publicly available preprint of Mati\'c, Radoi\v{c}i\'c, and Spomer
\cite{MRS} gives a computer-assisted proof of the upper bound seven.
Only the explicitly identified specialization in
Corollary~\ref{cor:codim-seven} uses that last result; our reductions and
product theorems are proved independently of it.

In arbitrary ambient dimension, Solymosi and Zahl proved
\begin{equation}\label{eq:SZ}
 \mathcal{N}(1,n)\leqslant4n\,7^{2n-3}\qquad(n\geqslant2)
\end{equation}
by expressing the contact relation through polynomial equalities and
non-equalities and using restricted Zariski closure \cite{SZ}.
General finiteness for arbitrary $k$ was established in the author's
preprint \cite{Li}, with the explicit estimate
$\mathcal{N}(k,n)\leqslant2^{\,2^{\,10(k+1)(n-k)}}$.
That argument uses sums of squares of minors, connected Nash pieces,
and descent through the possible direction-sum ranks.

The present paper combines a global-coordinate reformulation of that
finiteness mechanism with explicit bounds and structural relations between
the parameters. The role of H\"ofer's method is specific: we extend its
polynomial linearization and total-isotropy viewpoint, not its special
three-dimensional independence lemma. The passage from isotropy to a
cardinality bound uses real algebraic geometry. The explicit decomposition
estimate is proved below, so that no unspecified stratification constant
remains in the final bound. The reduction and product results supply,
respectively, exact relations and constructive lower bounds.

\subsection{Explicit upper bounds}
Set
\begin{equation}\label{eq:m-h}
 m=(k+1)(n-k),\qquad h=\min\{k,n-k-1\}.
\end{equation}
The number $m$ is the dimension of an affine Grassmann chart, whereas $h$
is the number of rank levels above the parallel rank.

\begin{theorem}[Explicit finiteness]\label{thm:main}
For every $0\leqslant k<n$,
\begin{equation}\label{eq:main-bound}
 \mathcal{N}(k,n)\leqslant\binom nk(n-k+1)
 (32n)^{\,h(4^m-1)/3}
 \leqslant 2^{\,2^{\,3m}}.
\end{equation}
In particular, every family of affine $k$-planes in $\R^n$ with a common
strictly positive separation distance is finite.
\end{theorem}

The estimate is deliberately nonoptimal. It improves the numerical
constant in the double-exponential bound of \cite{Li}, but does not change
its double-exponential nature. For lines, \eqref{eq:SZ} is much stronger.
The global projection model is used to explain the argument; a separate
chart calculation is necessary to obtain the exponent in
\eqref{eq:main-bound}.

\subsection{Exact reductions}
For $k\leqslant r\leqslant n$, write $\mathcal{N}_{\leqslant r}(k,n)$ for the extremal number in
\eqref{eq:N-definition} with the additional condition
\begin{equation}\label{eq:restricted-rank}
 \dim(U_L+U_{L'})\leqslant r\qquad(L,L'\in\cL),
\end{equation}
where $U_L$ is the direction of $L$.

\begin{theorem}[Rank-one excess and codimension-two reduction]
\label{thm:reductions}
For $1\leqslant k<n$,
\begin{equation}\label{eq:rank-one-reduction}
 \mathcal{N}_{\leqslant k+1}(k,n)=\mathcal{N}(1,n-k+1).
\end{equation}
Moreover, if $n\geqslant3$ and $\cL$ is an equidistant family of affine
$(n-2)$-planes in $\R^n$, there is an orthogonal decomposition
$\R^n=E\oplus K$, with $\dim E=3$ and $\dim K=n-3$, such that
\begin{equation}\label{eq:codim-splitting}
 L=\ell_L+K\quad(L\in\cL),
 \qquad d(L,L')=d(\ell_L,\ell_{L'}),
\end{equation}
where each $\ell_L$ is an affine line in $E$. Consequently,
\begin{equation}\label{eq:codim-equality}
 \mathcal{N}(n-2,n)=\mathcal{N}(1,3)\qquad(n\geqslant3).
\end{equation}
\end{theorem}

The direction-space dichotomy underlying the proof is the classical
star--top phenomenon: a family of $k$-subspaces meeting pairwise in
codimension at most one either contains a common $(k-1)$-subspace or lies
in a common $(k+1)$-subspace. An explicit finite-field formulation appears
in \cite[Proposition~17 and Corollary~18]{GR}; its dimension argument is
field-independent. We give an elementary proof over $\R$ and then establish
the additional Euclidean assertions. The finite-field subspace distance
in \cite{GR} is different from \eqref{eq:distance-definition}.
We also prove the sharp bound $2(n-k)$ in the top case
(Theorem~\ref{thm:top}).

\subsection{Multiplicative constructions}
\begin{theorem}[Transverse product]\label{thm:product}
Let $0\leqslant k_i<n_i$ for $i=1,2$. If $2k_1\geqslant n_1$, then
\begin{equation}\label{eq:product-main}
 \mathcal{N}(k_1+k_2,n_1+n_2)\geqslant \mathcal{N}(k_1,n_1)\mathcal{N}(k_2,n_2).
\end{equation}
The same conclusion holds if $2k_2\geqslant n_2$.
\end{theorem}

The construction rotates the second factor independently over each member
of the first family. Distances within one group come from the second
factor; distances between groups come from the first factor. Thus it is
not the ordinary Cartesian product of two equidistant configurations.
The simplex-and-line-pair construction in \cite[Section~3]{DKL} is a
special instance of the same allocation of distances.

Put
\begin{equation}\label{eq:a-S-intro}
 a=\mathcal{N}(1,3),\qquad S(c)=\mathcal{N}(c-1,2c-1)\quad(c\geqslant1).
\end{equation}
The product theorem implies
\begin{equation}\label{eq:super-intro}
 S(c+d)\geqslant S(c)S(d),\qquad
 S(c)\geqslant2^{c\bmod2}a^{\lfloor c/2\rfloor}.
\end{equation}
The bound $a\geqslant7$ from \cite{BLR} therefore yields, for example,
$\mathcal{N}(2,5)\geqslant14$, $\mathcal{N}(3,7)\geqslant49$, and $\mathcal{N}(5,11)\geqslant343$.
These lower bounds do not use the external upper bound $a\leqslant7$.

\subsection{Fixed-codimension finiteness and stability}

The preceding upper bound is uniform only after both $k$ and $n$ are
specified. A further consequence of the normal-space geometry is that,
for each fixed codimension, the extremal numbers are uniformly bounded
as the ambient dimension varies. For $c\geqslant1$, define
\[
 b_c:=\sup_{k\geqslant0}\mathcal{N}(k,k+c),
\]
where the supremum is initially allowed to be infinite.

\begin{theorem}[Fixed-codimension finiteness and eventual stability]
\label{thm:intro-stability}
For every integer $c\geqslant1$, one has $b_c<\infty$. Moreover, there
exists an integer $D_c\geqslant c$ such that
\[
 \mathcal{N}(k,k+c)\leqslant\mathcal{N}(D_c-c,D_c)
 \qquad(k\geqslant0),
\]
and
\[
 \boxed{\;
 \mathcal{N}(k,k+c)=\mathcal{N}(D_c-c,D_c)
 \qquad(k\geqslant D_c-c).
 \;}
\]
One may choose $D_1=1$ and, inductively,
\begin{equation}\label{eq:intro-Dc}
 D_c=c+\sum_{t=1}^{c-1}t(c-t)b_t
 \qquad(c\geqslant2).
\end{equation}
\end{theorem}

Thus, for each fixed codimension, the nondecreasing integer sequence
$k\mapsto\mathcal{N}(k,k+c)$ is eventually constant. The proof is based
on a distance-preserving projection lemma for normal spaces: a finite
configuration of codimension $c$ can be represented, without changing
its cardinality or pairwise distances, in an ambient space whose
dimension is bounded in terms of the lower-codimension extremal numbers.

In codimension three, the general recursive dimension bound can be
improved substantially. Put $a=\mathcal{N}(1,3)$.

\begin{theorem}[Codimension-three stability]
\label{thm:intro-codim-three}
Every finite unit-distance family of affine subspaces of codimension
three in any Euclidean space satisfies
\[
 \dim\left(\sum_i U_i^\perp\right)\leqslant a+3.
\]
Consequently,
\begin{equation}\label{eq:intro-codim-three-stability}
 \mathcal{N}(k,k+3)\leqslant\mathcal{N}(a,a+3)
 \qquad(k\geqslant0),
\end{equation}
and equality holds whenever $k\geqslant a$:
\[
 \boxed{\;
 \mathcal{N}(k,k+3)=\mathcal{N}(a,a+3)
 \qquad(k\geqslant a).
 \;}
\]
Using the computer-assisted result $\mathcal{N}(1,3)=7$, this gives
\[
 \boxed{\;
 \mathcal{N}(k,k+3)=\mathcal{N}(7,10)
 \qquad(k\geqslant7).
 \;}
\]
The bound $a+3$ for the span of the normal spaces is sharp.
\end{theorem}

Finally, the stabilized values themselves inherit the transverse
product inequality. In particular, the sequence $(b_c)_{c\geqslant1}$
is supermultiplicative:
\begin{equation}\label{eq:intro-bc-super}
 \boxed{\quad b_{c+d}\geqslant b_cb_d
 \qquad(c,d\geqslant1).\quad}
\end{equation}
Since $b_1=2$ and $b_2=a$, it follows that
$b_c\geqslant2^{c\bmod2}a^{\lfloor c/2\rfloor}$.

\subsection{Organization}
Sections~\ref{sec:projection}-\ref{sec:veronese} develop projection
coordinates, distance kernels, and their Veronese linearization.
Section~\ref{sec:nash} proves an explicit Nash decomposition bound.
Sections~\ref{sec:descent}-\ref{sec:charts} establish finiteness and
\eqref{eq:main-bound}. Sections~\ref{sec:basic}-\ref{sec:reductions}
prove the structural equalities and the sharp top bound.
Sections~\ref{sec:products}-\ref{sec:consequences} develop the product
inequalities and their consequences. Section~\ref{sec:stability} proves
the fixed-codimension and codimension-three stability theorems stated
above. The final discussion records the remaining quantitative and
structural questions.

\section{A global orthogonal-projection model}\label{sec:projection}

Let $L=b+U$, and let $P$ be the orthogonal projection onto $U$. Set
$a=(I_n-P)b$. Then $L=a+U$ and $a\in U^\perp$. The vector $a$ is independent
of the choice of $b\in L$, and it is the unique point of $L$ nearest to the
origin, since
\[
 \norm{a+u}^2=\norm{a}^2+\norm{u}^2\qquad(u\in U).
\]
We call $a$ the \emph{canonical point} of $L$.

\begin{proposition}\label{prop:model}
The affine Grassmannian $\Graff(k,n)$ is identified with the real algebraic
set
\begin{equation}\label{eq:model}
 \A_{k,n}=\left\{(P,a)\in\Sym_n(\R)\times\R^n\Big|
 P^2=P,\ \tr P=k,\ Pa=0\right\}
\end{equation}
by $(P,a)\mapsto a+\im P$. The ambient coordinate space has dimension
$M_n=n(n+3)/2$.
\end{proposition}

\begin{proof}
An orthogonal projection onto a $k$-dimensional subspace is symmetric,
idempotent and has trace $k$. Its kernel is the orthogonal complement of
its image, so $Pa=0$ is equivalent to $a\perp\im P$.
Conversely, a real symmetric idempotent matrix has only the eigenvalues zero
and one. Thus the trace condition makes its rank equal to $k$, and it is the
orthogonal projection onto its image. The uniqueness of the direction and
of the canonical point proves bijectivity. Finally,
$\dim\Sym_n(\R)+n=n(n+1)/2+n=M_n$.
\end{proof}

Only the independent upper-triangular entries of $P$ are used as coordinates,
so symmetry is built into the ambient space. The equations in
\eqref{eq:model} have degree at most two. In particular, no inequalities are
needed to impose positive semidefiniteness: it follows from symmetry and
idempotence over $\R$.

For completeness, this model agrees with the usual Nash structure of the
affine Grassmannian. Fix a coordinate splitting
$\R^n=\R^k\oplus\R^{n-k}$. On a graph chart, put
\[
 U_B=\left\{(u,Bu)\Big|u\in\R^k\right\},\qquad B\in\R^{(n-k)\times k}.
\]
For $c\in\R^{n-k}$, the corresponding projection and canonical point are
\begin{equation}\label{eq:graph-chart}
 \begin{split}
 P_B&=\begin{pmatrix}I_k\\B\end{pmatrix}
 (I_k+B\T B)^{-1}\begin{pmatrix}I_k&B\T\end{pmatrix},\\
 a_{B,c}&=\begin{pmatrix}-B\T c\\c\end{pmatrix}.
 \end{split}
\end{equation}
The inverse matrix is a Nash function of $B$, because $I_k+B\T B$ is
positive definite. Conversely, if $P$ is written in blocks, then
$B=P_{21}P_{11}^{-1}$, and $c$ is the last $n-k$ coordinates of $a$.
Permuting the coordinate splitting gives a finite collection of graph
charts covering the model. Hence $\A_{k,n}$ is a smooth Nash manifold of
intrinsic dimension
\begin{equation}\label{eq:intrinsic-dimension}
 (k+1)(n-k).
\end{equation}
The endpoint $k=0$ is interpreted directly as $P=0$, $a\in\R^n$.
The intrinsic dimension \eqref{eq:intrinsic-dimension} should be distinguished
from the ambient dimension $M_n$ used in our decomposition bound.

The canonical-point coordinate gives the elementary estimate
\begin{equation}\label{eq:canonical-distance-bound}
 d(L_{(P,a)},L_{(Q,b)})\leqslant\norm{a-b}.
\end{equation}
This estimate, rather than continuity of separation distance, is the
geometric input in the diagonal argument of Section~\ref{sec:descent}.

\section{A distance identity on each rank stratum}\label{sec:distance}

Let $x=(P,a)$ and $y=(Q,b)$ belong to $\A_{k,n}$, and write
\begin{equation}\label{eq:S-w}
 S=P+Q,\qquad w=b-a.
\end{equation}
For every $v\in\R^n$,
\[
 v\T Sv=\norm{Pv}^2+\norm{Qv}^2.
\]
It follows that $S\succeq0$ and
\begin{equation}\label{eq:sum-directions}
 \ker S=\ker P\cap\ker Q,
 \qquad \im S=\im P+\im Q.
\end{equation}
The difference set of the two affine planes is
$w+(\im P+\im Q)$. Consequently,
\begin{equation}\label{eq:projection-distance}
 d(L_x,L_y)^2
 =\dist(w,\im S)^2
 =\norm{\proj_{\ker S}w}^2.
\end{equation}
In particular, positive separation implies $\rank S\leqslant n-1$.

For $A\in\Sym_n(\R)$, define
\begin{equation}\label{eq:elementary}
 e_0(A)=1,\qquad
 e_j(A)=\sum_{\substack{I\subset\{1,\ldots,n\}\\|I|=j}}
                  \det A_{I,I}\quad(1\leqslant j\leqslant n).
\end{equation}
Equivalently,
\[
 \det(tI_n+A)=\sum_{j=0}^n e_j(A)t^{n-j}.
\]
Thus $e_j(A)$ is the $j$th elementary symmetric polynomial of the
eigenvalues and is invariant under orthogonal changes of coordinates.

\begin{lemma}\label{lem:rank-test}
If $A\succeq0$, then for $1\leqslant j\leqslant n$,
\begin{equation}\label{eq:rank-test}
 e_j(A)>0\iff\rank A\geqslant j,
 \qquad e_j(A)=0\iff\rank A<j.
\end{equation}
\end{lemma}

\begin{proof}
All eigenvalues of $A$ are nonnegative. Each term in their $j$th elementary
symmetric polynomial is nonnegative, and a positive term occurs precisely
when at least $j$ eigenvalues are positive.
\end{proof}

\begin{lemma}[Rank-one update identity]\label{lem:rank-one-update}
Let $S\succeq0$ have rank $r<n$, and let $w\in\R^n$. Then
\begin{equation}\label{eq:rank-one-update}
 e_{r+1}(S+ww\T)=e_r(S)\norm{\proj_{\ker S}w}^2.
\end{equation}
\end{lemma}

\begin{proof}
If $r=0$, then $S=0$ and both sides equal $\norm{w}^2$. Suppose $r\geqslant1$.
By an orthogonal change of coordinates, we may assume
\[
 S=\diag(\lambda_1,\ldots,\lambda_r,0,\ldots,0),\qquad \lambda_i>0.
\]
Then $e_r(S)=\lambda_1\cdots\lambda_r$.

Consider a principal minor of order $r+1$ of $S+ww\T$. If its index set
does not contain all of $1,\ldots,r$, the corresponding submatrix of $S$
has rank at most $r-1$. Adding a rank-one matrix raises the rank to at most
$r$, so that minor vanishes. The remaining index sets are exactly
$\{1,\ldots,r,j\}$ with $j>r$. For each such set, the determinant is
\[
 \lambda_1\cdots\lambda_r\,w_j^2.
\]
For instance, this follows from
$\det(A+uu\T)=\det A+u\T\adj(A)u$, applied to
$A=\diag(\lambda_1,\ldots,\lambda_r,0)$.
Summing over $j>r$ proves \eqref{eq:rank-one-update}.
\end{proof}

\begin{proposition}[Global distance formula]\label{prop:distance-formula}
For $x=(P,a)$, $y=(Q,b)\in\A_{k,n}$ with $r=\rank(P+Q)<n$,
\begin{equation}\label{eq:distance-formula}
 d(L_x,L_y)^2
 =\frac{e_{r+1}\bigl(P+Q+(b-a)(b-a)\T\bigr)}{e_r(P+Q)}.
\end{equation}
The denominator is positive on the rank-$r$ stratum, with the convention
$e_0=1$.
\end{proposition}

\begin{proof}
Combine \eqref{eq:projection-distance} with
Lemma~\ref{lem:rank-one-update} and Lemma~\ref{lem:rank-test}.
\end{proof}

\begin{remark}\label{rem:rank-dependence}
Formula~\eqref{eq:distance-formula} does not use a choice of basis for either
direction or a choice of nonvanishing minor. Its denominator nevertheless
depends on the rank. On a lower-rank stratum, both the numerator and the
denominator of a higher-rank formula can vanish. Such vanishing does not
encode a distance value. The kernels below retain this distinction.
\end{remark}

\section{Polynomial kernels adapted to rank degeneration}\label{sec:kernels}

Throughout the rank-descent argument, set
\begin{equation}\label{eq:rstar}
 r_* = \min\{2k,n-1\}.
\end{equation}
If distinct planes in a family have positive pairwise separation, then
\eqref{eq:sum-directions} and \eqref{eq:projection-distance} imply
\[
 k\leqslant\rank(P_x+P_y)\leqslant r_*.
\]
If $r_*=k$, no descent is needed. Otherwise, fix an integer
$k<\rho\leqslant r_*$ and define
\begin{align}
 F_\rho(x,y)
   &=e_{\rho+1}\bigl(P+Q+(b-a)(b-a)\T\bigr)-e_\rho(P+Q),
      \label{eq:F-kernel}\\
 R_\rho(x,y)&=e_{\rho+1}(P+Q).
      \label{eq:R-kernel}
\end{align}
These are polynomials on the whole ambient coordinate space
$\R^{M_n}\times\R^{M_n}$, not only on $\A_{k,n}\times\A_{k,n}$.
They are symmetric in $x$ and $y$.

\begin{proposition}\label{prop:encoding}
For $x,y\in\A_{k,n}$, the following assertions hold.
\begin{enumerate}
\item If $\rank(P_x+P_y)\leqslant\rho$ and $d(L_x,L_y)=1$, then
$F_\rho(x,y)=R_\rho(x,y)=0$.
\item If $\rank(P_x+P_y)\leqslant\rho-1$, both kernels vanish regardless of the
distance. In particular, $F_\rho(x,x)=R_\rho(x,x)=0$ for every $x$.
\item If $F_\rho(x,y)=R_\rho(x,y)=0$ and
$e_\rho(P_x+P_y)>0$, then
$\rank(P_x+P_y)=\rho$ and $d(L_x,L_y)=1$.
\end{enumerate}
\end{proposition}

\begin{proof}
Write $S=P_x+P_y$ and $w=a_y-a_x$.
If $\rank S=\rho$ and the distance is one, then
Proposition~\ref{prop:distance-formula} gives $F_\rho=0$, while
Lemma~\ref{lem:rank-test} gives $R_\rho=0$.
If $\rank S\leqslant\rho-1$, then
\[
 e_\rho(S)=e_{\rho+1}(S)=0,
 \qquad \rank(S+ww\T)\leqslant\rho,
\]
so $e_{\rho+1}(S+ww\T)=0$ as well. This proves the first two assertions,
since on the diagonal $S=2P_x$ has rank $k<\rho$ and $w=0$.

For the third assertion, the vanishing of
$R_\rho=e_{\rho+1}(S)$ implies $\rank S\leqslant\rho$, whereas
$e_\rho(S)>0$ implies $\rank S\geqslant\rho$. Thus $\rank S=\rho$.
The identity $F_\rho=0$ and Proposition~\ref{prop:distance-formula} now give
$d(L_x,L_y)^2=1$.
\end{proof}

\begin{lemma}[Degree bounds]\label{lem:kernel-degrees}
The total degrees in both groups of ambient coordinates satisfy
\begin{equation}\label{eq:kernel-degrees}
 \deg R_\rho\leqslant\rho+1,
 \qquad \deg F_\rho\leqslant\rho+2\leqslant n+1.
\end{equation}
\end{lemma}

\begin{proof}
The polynomial $e_{\rho+1}(S)$ is homogeneous of degree $\rho+1$ in the
entries of $S$, and $S=P+Q$ is linear in the ambient coordinates.
For each principal index set $I$ of size $\rho+1$, the polynomial identity
\[
 \det(S_{I,I}+w_Iw_I\T)
 =\det S_{I,I}+w_I\T\adj(S_{I,I})w_I
\]
holds whether or not $S_{I,I}$ is invertible. The adjugate entries have
degree $\rho$ in $S$, and the update term is quadratic in $w$. Since $w=b-a$
is linear, this term has total degree at most $\rho+2$. Summing the
principal minors and subtracting $e_\rho(S)$ proves the claim.
\end{proof}

\begin{remark}[The role of rank one]\label{rem:rank-one-degree}
The update $ww\T$ has rank at most one in every ambient dimension, since
$w=b-a$ is a single vector. In the identity
\[
 e_{\rho+1}(S+ww\T)=e_{\rho+1}(S)
 +\sum_{|I|=\rho+1}w_I\T\adj(S_{I,I})w_I,
\]
the first term has $S$-degree $\rho+1$; only the coefficients of the
quadratic $w$-terms have $S$-degree $\rho$.
This distinction matters. If one replaces $ww\T$ by $WW\T$ with $r$ columns,
multilinearity in the columns gives only
\[
 \deg_{S,W} e_j(S+WW\T)\leqslant j+\min\{r,j\}.
\]
For example, $S=\diag(s_1,s_2)$ and $W=\diag(u,v)$ produce the fourth-degree
term $u^2v^2$ in $\det(S+WW\T)$. Such a replacement is not made here.
\end{remark}

The upper-rank kernel $R_\rho$ is essential in the recovery direction:
$F_\rho=0$ alone would not justify applying the rank-$\rho$ distance
formula. Retaining the two kernels separately also avoids squaring the
unit-distance kernel before the Veronese lifting.

\section{Veronese lifting and algebraic linear sections}\label{sec:veronese}

Put
\begin{equation}\label{eq:veronese-dimension}
 D_n=\binom{M_n+n+1}{n+1},
\end{equation}
and let
\begin{equation}\label{eq:veronese}
 \Phi_{n+1}:\R^{M_n}\longrightarrow\R^{D_n},\qquad
 \Phi_{n+1}(x)=(x^\alpha)_{|\alpha|\leqslant n+1}.
\end{equation}
We include the constant monomial. Thus this is the affine, nonhomogeneous
version of the Veronese map.

\begin{lemma}[Bilinearization]\label{lem:bilinearization}
Suppose $K(x,y)=K(y,x)$ is a real polynomial whose degree in each group of
variables is at most $n+1$. Then there is a symmetric bilinear form
$\B_K$ on $\R^{D_n}$ such that
\begin{equation}\label{eq:bilinearization}
 K(x,y)=\B_K\bigl(\Phi_{n+1}(x),\Phi_{n+1}(y)\bigr).
\end{equation}
\end{lemma}

\begin{proof}
Expand uniquely
\[
 K(x,y)=\sum_{|\alpha|,|\beta|\leqslant n+1}
                  c_{\alpha\beta}x^\alpha y^\beta.
\]
Symmetry gives $c_{\alpha\beta}=c_{\beta\alpha}$. The coefficient matrix
$C_K=(c_{\alpha\beta})$ represents the desired form:
$\B_K(u,v)=u\T C_Kv$.
\end{proof}

By Lemma~\ref{lem:kernel-degrees}, this construction applies simultaneously
to $F_\rho$ and $R_\rho$. Write the resulting forms as
$\B_\rho^F$ and $\B_\rho^R$. They need not be nondegenerate. We will use
neither their signatures nor their Witt indices.

Let $X\subset\A_{k,n}$ be a finite unit-distance configuration with
$\rank(P_x+P_y)\leqslant\rho$ for all $x,y\in X$, and set
\begin{equation}\label{eq:lifted-span}
 W_X=\Span\left\{\Phi_{n+1}(x)\big|x\in X\right\}.
\end{equation}
Proposition~\ref{prop:encoding} includes the diagonal and therefore gives
\begin{equation}\label{eq:common-isotropy}
 \B_\rho^F\big|_{W_X\times W_X}=0,
 \qquad
 \B_\rho^R\big|_{W_X\times W_X}=0.
\end{equation}
In other words, $W_X$ is simultaneously totally isotropic for the two
forms. Define the associated algebraic linear section by
\begin{equation}\label{eq:linear-section}
 Z_X=\left\{x\in\A_{k,n}\big|\Phi_{n+1}(x)\in W_X\right\}.
\end{equation}
Here ``linear section'' refers to a linear section in the lifted space,
followed by inverse image in the parameter space. We have $X\subset Z_X$.

\begin{proposition}[Propagation of the polynomial relations]\label{prop:propagation}
For every $x,y\in Z_X$,
\begin{equation}\label{eq:propagation}
 F_\rho(x,y)=R_\rho(x,y)=0.
\end{equation}
\end{proposition}

\begin{proof}
Both lifted points lie in $W_X$, so the assertion follows directly from
\eqref{eq:common-isotropy} and \eqref{eq:bilinearization}.
\end{proof}

Equivalently, $Z_X$ is cut out inside $\A_{k,n}$ by all polynomials of degree
at most $n+1$ that vanish on $X$. Indeed, the coefficient vectors of
these polynomials form the linear annihilator of $W_X$. This observation
also makes the bounded complexity of $Z_X$ transparent.
The propagated identities are polynomial identities, not unconditional
unit-distance identities: the latter are recovered only where
$e_\rho(P_x+P_y)>0$.

\begin{proposition}[Uniform algebraic complexity]\label{prop:complexity}
The set $Z_X$ is the real zero set of one polynomial in $M_n$ variables of
degree at most $2n+2$. The degree bound is independent of the cardinality
and location of $X$.
\end{proposition}

\begin{proof}
Choose a basis $\lambda_1,\ldots,\lambda_t$ of the annihilator of $W_X$,
where $0\leqslant t\leqslant D_n$, and put
$f_j(x)=\lambda_j(\Phi_{n+1}(x))$. Then $\deg f_j\leqslant n+1$, and
$\Phi_{n+1}(x)\in W_X$ is equivalent to $f_j(x)=0$ for every $j$.
Define
\begin{equation}\label{eq:single-polynomial}
 \begin{split}
 \cP_X(P,a)={}&\sum_{j=1}^t f_j(P,a)^2
       +\norm{P^2-P}_F^2\\
       &+(\tr P-k)^2+\norm{Pa}^2.
 \end{split}
\end{equation}
Over $\R$, a sum of squares vanishes precisely when each summand vanishes.
Consequently, \eqref{eq:model} and \eqref{eq:linear-section} give
$Z_X=\{\cP_X=0\}$. The last three terms have degree at most four, and the
first has degree at most $2n+2\geqslant4$.
\end{proof}

\begin{remark}\label{rem:sum-of-squares}
The sum-of-squares representation is used only to bound the complexity of
a real zero set. It does not make $Z_X$ a nonsingular hypersurface.
Indeed, the gradient of \eqref{eq:single-polynomial} vanishes throughout its
zero set. The next section applies to arbitrary, possibly singular, real
algebraic sets.
\end{remark}

\section{An explicit connected Nash decomposition bound}\label{sec:nash}

A \emph{Nash submanifold} is a semialgebraic embedded real-analytic
submanifold; it need not be closed. A Nash function is real analytic with
semialgebraic graph. Zero-dimensional pieces are allowed. We use the
Tarski--Seidenberg theorem and elementary facts about semialgebraic sets
as in \cite{BCR}. Rather than invoke an unspecified effective
stratification constant, we give a projection-and-lifting construction
with explicit estimates. No Whitney or frontier condition is needed.

\begin{lemma}[Root functions at constant multiplicity]\label{lem:roots}
Let $C$ be a connected Nash submanifold and let
$p(x,t)=\sum\limits_{j=0}^{a}p_j(x)t^j$ have Nash coefficients on $C$, with
$p_a$ nowhere zero. Suppose that the number of distinct complex roots of
$p(x,\cdot)$ is independent of $x\in C$. Then the number of distinct real
roots is constant, and the ordered real roots are Nash functions on $C$.
Their multiplicities are locally constant, and hence constant on $C$.
If a finite family of such polynomials and all its pairwise products have
constant numbers of distinct complex roots, any two of their real root
functions either coincide everywhere or are everywhere distinct.
\end{lemma}

\begin{proof}
Fix $x_0\in C$. Take pairwise disjoint small disks around the distinct
complex roots at $x_0$, with the disks at nonreal roots disjoint from the
real axis. Continuity of polynomial roots, with multiplicity, implies that
for nearby parameters all roots remain in these disks and the total
multiplicity in each disk is unchanged. Each disk contains at least one
distinct root. Constancy of the total number of distinct roots forces
exactly one distinct root in each disk. Thus no root splits and its
multiplicity is locally constant. A real root remains real: otherwise its
complex conjugate would be a second distinct root in the same disk.
Consequently, the number of real roots is locally constant and therefore
constant on connected $C$.

If a real root has multiplicity $\mu$, it is a simple root of
$\partial_t^{\mu-1}p$. The analytic implicit function theorem on $C$
therefore makes its continuation analytic. Ordering the real roots gives
globally defined functions whose local expressions agree. Their graphs
are semialgebraic: the $j$th root is characterized by the polynomial
equation and the existence of exactly $j-1$ distinct smaller real roots,
a first-order real-algebraic condition. Tarski--Seidenberg applies.
Thus the functions are Nash.

Finally, if roots of two polynomials coincide at $x_0$, apply the same
no-splitting argument to their product. They must coincide near $x_0$.
The coincidence set is therefore open and, by continuity, closed in $C$.
Connectedness proves the last assertion. Distinct real root functions
cannot cross without coinciding.
\end{proof}

\begin{theorem}[Quantitative Nash decomposition]\label{thm:quantitative-Nash}
Let $M\geqslant0$, $s,d\geqslant1$, and let $f_1,\ldots,f_s\in
\R[x_1,\ldots,x_M]$ have total degree at most $d$. There is a partition of
$\R^M$ into connected Nash submanifolds on each of which every $f_i$ has
constant sign. The number of pieces is at most
\begin{equation}\label{eq:CAD-bound}
 T(M,s,d)=(8sd)^{(4^M-1)/3}.
\end{equation}
Identically zero polynomials are permitted.
\end{theorem}

\begin{proof}
We induct on $M$. For $M=0$ there is a single point and the exponent is
zero. Suppose $M\geqslant1$, write $x=(y,t)$, and form the auxiliary family
\[
 \mathcal G=\{f_i|1\leqslant i\leqslant s\}
 \cup\{f_if_j|1\leqslant i<j\leqslant s\}.
\]
It has at most $s^2$ members, each of degree at most $2d$. Write
$g(y,t)=\sum\limits_{j=0}^{2d}g_j(y)t^j$ for each $g\in\mathcal G$.
Our projection family in the $M-1$ variables $y$ consists first of all
coefficients $g_j$. For every $1\leqslant a\leqslant2d$, also form the truncation
\[
 g^{[a]}(y,t)=\sum_{j=0}^{a}g_j(y)t^j
\]
and the Sylvester matrix
\[
 B_{g,a}(y)=\Syl\bigl(g^{[a]},\partial_t g^{[a]}\bigr),
\]
constructed with formal degrees $a$ and $a-1$. It is a
$(2a-1)\times(2a-1)$ matrix. Include in the projection family the polynomials
\begin{equation}\label{eq:projection-rank-polys}
 e_j\bigl(B_{g,a}B_{g,a}\T\bigr),\qquad 1\leqslant j\leqslant2a-1.
\end{equation}
The matrices are defined by their coefficient arrays even where a leading
coefficient vanishes. Only the truncation at the actual degree will be
used on a given base piece.

For each $g$ the number of projection polynomials is at most
\[
 (2d+1)+\sum_{a=1}^{2d}(2a-1)=4d^2+2d+1\leqslant7d^2.
\]
Each entry of $B_{g,a}$ has degree at most $2d$ in $y$, and an expression
in \eqref{eq:projection-rank-polys} has degree at most
$4dj\leqslant4d(4d-1)<16d^2$. Thus the projection family has size and degree
bounded by
\begin{equation}\label{eq:projection-cost}
 s'\leqslant7s^2d^2,\qquad d'\leqslant16d^2.
\end{equation}
Zero or repeated polynomials may be retained. If necessary, a constant
polynomial can be used when the nonzero projection family is empty.

Apply the induction hypothesis to the projection family. Let $C$ be a
resulting connected Nash piece. Since all coefficients $g_j$ have constant
sign on $C$, each specialization $g(y,\cdot)$ is either identically zero
for every $y\in C$, or has a fixed degree $a$ with nowhere-zero leading
coefficient. In the latter case $g=g^{[a]}$ on $C$. The signs of
\eqref{eq:projection-rank-polys} and Lemma~\ref{lem:rank-test} make the rank
of $B_{g,a}$ constant. The standard Sylvester identity, valid for a
polynomial of exact degree $a$ in characteristic zero, is
\begin{equation}\label{eq:sylvester-rank}
 \rank\Syl(g,g')=2a-1-\deg\gcd(g,g').
\end{equation}
Since $\deg g-\deg\gcd(g,g')$ is the number of its distinct complex roots,
that number is constant on $C$. A nonzero constant has no roots and
requires no Sylvester matrix.

Lemma~\ref{lem:roots}, applied also to the pairwise products, gives a
single ordered list of all distinct real roots of the nonzero $f_i$:
\[
 \theta_1(y)<\cdots<\theta_q(y),\qquad q\leqslant sd.
\]
All these functions are Nash, and each is either identically a root of a
given $f_i$ or never a root of it. Their graphs and the open bands between
successive graphs, including the two unbounded bands, partition
$C\times\R$ into at most $2sd+1$ connected Nash submanifolds. Graphs are
Nash by construction; bands are open subsets of $C\times\R$ and are
connected, being homeomorphic to $C\times(0,1)$. Every $f_i$ is either zero
or nowhere zero on each piece, so connectedness gives constant sign.
These pieces over all base pieces form the required partition.

It remains to count. Put $B=8sd$. By \eqref{eq:projection-cost},
\[
 8s'd'\leqslant896s^2d^4\leqslant(8sd)^4=B^4,
 \qquad 2sd+1\leqslant B.
\]
Hence the total number of lifted pieces is at most
\[
 B\,(B^4)^{(4^{M-1}-1)/3}
 =B^{(4^M-1)/3},
\]
which proves the induction.
\end{proof}

\begin{corollary}[One equation]\label{lem:uniform-nash}
The real zero set of any polynomial of degree at most $D\geqslant1$ in
$M\geqslant1$ variables is a disjoint union of at most
\begin{equation}\label{eq:Sigma}
 \Sigma(M,D)=(8D)^{(4^M-1)/3}
\end{equation}
connected Nash submanifolds. This includes the zero polynomial and the
empty zero set.
\end{corollary}
\begin{proof}
Apply Theorem~\ref{thm:quantitative-Nash} with $s=1$ and retain precisely
the zero-sign pieces. The exponent is an integer.
\end{proof}

\begin{remark}\label{rem:why-Nash}
The construction is an explicit version of cylindrical algebraic
decomposition, with a redundant projection family chosen for transparent
counting. We do not claim optimal complexity. The pieces need only be
connected and real analytic; no Whitney conditions are imposed.
Connected components of a singular algebraic set would not suffice,
because a polynomial can vanish on an open part of one branch but not on
another branch of the same connected component.
\end{remark}

\begin{lemma}[Real-analytic identity principle]\label{lem:analytic-identity}
Let $S$ be a connected real-analytic manifold and let $g:S\to\R$ be
analytic. If $g$ is not identically zero, its zero set has empty interior.
If also $g\geqslant0$, then $\{s\in S|g(s)>0\}$ is dense in $S$.
\end{lemma}
\begin{proof}
If $g$ vanishes on a nonempty open set, let $E$ be the set of points near
which it vanishes identically. This set is nonempty and open. It is closed:
at a limit point, all derivatives in a fixed analytic chart vanish by
continuity, and the convergent Taylor series makes $g$ zero locally.
Connectedness gives $E=S$. The assertions follow by contraposition.
\end{proof}

For the global model, we may now choose the explicit constant
\begin{equation}\label{eq:global-constant}
 M_n=\frac{n(n+3)}2,\qquad
 C_n=\Sigma(M_n,2n+2)
 =[16(n+1)]^{(4^{M_n}-1)/3}.
\end{equation}

\section{Rank descent on the linear section}\label{sec:descent}

Fix $k<\rho\leqslant r_*$ and a finite unit-distance configuration
$X\subset\A_{k,n}$ with direction-sum rank at most $\rho$.
By Proposition~\ref{prop:complexity} and Lemma~\ref{lem:uniform-nash}, there
is a partition
\begin{equation}\label{eq:Z-partition}
 Z_X=S_1\sqcup\cdots\sqcup S_s,
 \qquad s\leqslant C_n,
\end{equation}
where each nonempty $S_j$ is a connected Nash submanifold.

\begin{theorem}[Strict descent on each Nash piece]\label{thm:rank-descent}
For every piece $S_j$ in \eqref{eq:Z-partition} and all $x,y\in S_j$,
\begin{equation}\label{eq:rank-descent}
 \rank(P_x+P_y)\leqslant\rho-1.
\end{equation}
\end{theorem}

\begin{proof}
Fix a nonempty piece $S$. Proposition~\ref{prop:propagation} gives
\[
 F_\rho(x,y)=0,
 \qquad e_{\rho+1}(P_x+P_y)=R_\rho(x,y)=0
 \qquad(x,y\in S).
\]
Since $P_x+P_y$ is positive semidefinite,
Lemma~\ref{lem:rank-test} yields
\begin{equation}\label{eq:rank-upper-on-S}
 \rank(P_x+P_y)\leqslant\rho\qquad(x,y\in S).
\end{equation}
Suppose that equality holds for some $x_0,y_0\in S$.
On the connected real-analytic manifold $S\times S$, consider
\begin{equation}\label{eq:analytic-rank-function}
 g(x,y)=e_\rho(P_x+P_y).
\end{equation}
This function is analytic and nonnegative, with $g(x_0,y_0)>0$.
Lemma~\ref{lem:analytic-identity} therefore makes $\{(x,y)\in S\times S|g(x,y)>0\}$ dense in
$S\times S$.

Choose $z\in S$. There are sequences $x_\nu,y_\nu\in S$ such that
\[
 (x_\nu,y_\nu)\longrightarrow(z,z),
 \qquad g(x_\nu,y_\nu)>0.
\]
Both kernels vanish at these pairs. The recovery assertion in
Proposition~\ref{prop:encoding} gives
\[
 \rank(P_{x_\nu}+P_{y_\nu})=\rho,
 \qquad d(L_{x_\nu},L_{y_\nu})=1.
\]
The pairs are distinct because the diagonal has rank $k<\rho$.
Let $a_{x_\nu},a_{y_\nu}$ denote their canonical points. Ambient coordinate
convergence and \eqref{eq:canonical-distance-bound} imply
\begin{equation}\label{eq:diagonal-contradiction}
 1=d(L_{x_\nu},L_{y_\nu})
 \leqslant\norm{a_{x_\nu}-a_{y_\nu}}\longrightarrow0,
\end{equation}
a contradiction. Thus equality never occurs in
\eqref{eq:rank-upper-on-S}, proving \eqref{eq:rank-descent}.
\end{proof}

The conclusion holds on the entire square of each Nash piece, not merely
on the original configuration points in that piece. It is this uniform
loss of one rank level that allows the argument to be iterated.

\begin{remark}\label{rem:continuity}
Separation distance on the affine Grassmannian is not generally continuous
at a rank drop. No such continuity is used in
\eqref{eq:diagonal-contradiction}. The right side is a continuous function
of the canonical-point coordinates and vanishes on the diagonal. The
left side is forced to be one only at the chosen full-rank pairs, by the
polynomial relations and their recovery condition.
\end{remark}

\section{The recursion and proof of finiteness}\label{sec:finiteness}

For $k\leqslant\rho\leqslant r_*$, let $A_\rho$ denote the supremum of $|X|$ over finite
sets $X\subset\A_{k,n}$ satisfying
\begin{equation}\label{eq:A-definition}
 \begin{cases}
 d(L_x,L_y)=1,&x\ne y,\\
 \rank(P_x+P_y)\leqslant\rho,&x,y\in X.
 \end{cases}
\end{equation}
We initially allow $A_\rho=+\infty$ and prove finiteness by induction.

\begin{lemma}[The parallel rank]\label{lem:parallel}
One has $A_k\leqslant n-k+1$.
\end{lemma}

\begin{proof}
If \eqref{eq:A-definition} holds with $\rho=k$, then the $k$-dimensional
spaces $\im P_x$ all coincide, because the sum of any two has dimension at
most $k$. Denote the common direction by $U$.
The canonical points $a_1,\ldots,a_m$ lie in $U^\perp\cong\R^{n-k}$, and
\[
 \norm{a_i-a_j}=d(L_i,L_j)=1\qquad(i\ne j).
\]
For $m\geqslant2$, put $u_i=a_i-a_m$, $1\leqslant i<m$. Then
$\norm{u_i}^2=1$ and $u_i\cdot u_j=1/2$ for $i\ne j$. Their Gram matrix is
\[
 \frac12(I_{m-1}+J_{m-1}),
\]
where $J_{m-1}$ is the all-ones matrix. This matrix is positive definite,
so the $u_i$ are linearly independent. Hence $m-1\leqslant n-k$.
The case $m=1$ is immediate.
\end{proof}

\begin{proposition}[Rank recursion]\label{prop:rank-recursion}
For $k<\rho\leqslant r_*$, if $A_{\rho-1}<\infty$, then
\begin{equation}\label{eq:rank-recursion}
 A_\rho\leqslant C_n A_{\rho-1}.
\end{equation}
\end{proposition}

\begin{proof}
Take any finite configuration $X$ counted by $A_\rho$ and construct
$Z_X$ and its partition \eqref{eq:Z-partition}.
By Theorem~\ref{thm:rank-descent}, each $X\cap S_j$ has direction-sum rank
at most $\rho-1$. It retains the unit-distance condition, so
$|X\cap S_j|\leqslant A_{\rho-1}$. Therefore
\[
 |X|=\sum_{j=1}^s|X\cap S_j|
 \leqslant sA_{\rho-1}\leqslant C_nA_{\rho-1}.
\]
Taking the supremum over $X$ proves the proposition.
\end{proof}

\begin{proposition}[Global-coordinate bound]\label{prop:global-bound}
With $C_n$ as in \eqref{eq:global-constant}, one has
\begin{equation}\label{eq:global-bound}
 \mathcal{N}(k,n)\leqslant(n-k+1)C_n^{\min\{k,n-k-1\}}.
\end{equation}
\end{proposition}
\begin{proof}
Starting with Lemma~\ref{lem:parallel} and applying
Proposition~\ref{prop:rank-recursion} successively yields
\[
 A_{r_*}\leqslant(n-k+1)C_n^{r_*-k}.
\]
Every unit-distance configuration has direction-sum rank at most
$r_*=\min\{2k,n-1\}$ for distinct pairs, by
\eqref{eq:sum-directions} and \eqref{eq:projection-distance}; its diagonal
has rank $k\leqslant r_*$. Hence it is counted by $A_{r_*}$, and
$\mathcal{N}(k,n)\leqslant A_{r_*}$. Since
\[
 r_*-k=\min\{k,n-k-1\},
\]
this is \eqref{eq:global-bound}. If $r_*=k$, the iteration is empty and the
parallel bound applies directly. This includes $k=0$ and $k=n-1$.
\end{proof}

\begin{corollary}\label{cor:endpoints}
For every $n\geqslant1$,
\[
 \mathcal{N}(0,n)=n+1,\qquad \mathcal{N}(n-1,n)=2.
\]
\end{corollary}

\begin{proof}
The upper bounds follow from the zero-exponent cases of
Proposition~\ref{prop:global-bound}. A regular simplex of side length one gives the first
lower bound, and two parallel hyperplanes at separation one give the second.
For $n=1$, both statements read $\mathcal{N}(0,1)=2$.
\end{proof}

\begin{remark}\label{rem:product}
The multiplicative recursion reflects the conclusion of rank descent:
a Nash piece can contain a whole lower-rank configuration, not just a
single original point. There are at most $C_n$ pieces, and each can contain
at most $A_{\rho-1}$ points. No bound on relations between different pieces
is needed.
\end{remark}

\section{Explicit bounds in affine Grassmann charts}\label{sec:charts}

The global proof uses $M_n=n(n+3)/2$ ambient coordinates. Replacing that
number by the intrinsic dimension without changing coordinates would be
unjustified. We now repeat the polynomial construction on genuine graph
charts and track the resulting degrees.

Fix a set of $k$ coordinate axes and permute coordinates so that they come
first. A plane in the associated chart is uniquely written as
\begin{equation}\label{eq:affine-graph}
 L_{(B,b)}=\left\{(u,Bu+b)\Big|u\in\R^k\right\},
 \quad B\in\R^{(n-k)\times k},\quad b\in\R^{n-k}.
\end{equation}
The parameter space is all of $\R^m$, where $m=(k+1)(n-k)$.
There are $\binom nk$ such charts, and they cover the affine Grassmannian,
since every rank-$k$ direction matrix has a nonzero $k\times k$ minor.
The point $p_{(B,b)}=(0,b)$ lies on the plane, so
\begin{equation}\label{eq:graph-point-distance}
 d(L_{(B,b)},L_{(C,c)})\leqslant\norm{b-c}.
\end{equation}

For $x=(B,b)$ and $y=(C,c)$ in a single chart, put
\begin{equation}\label{eq:chart-Gram}
 A(x,y)=\begin{pmatrix}I_k&I_k\\B&C\end{pmatrix},
 \quad G(x,y)=A(x,y)A(x,y)\T,
 \quad w(x,y)=\begin{pmatrix}0\\c-b\end{pmatrix}.
\end{equation}
Then $G\succeq0$, its image is $U_x+U_y$, and
$\rank G=\dim(U_x+U_y)$. The entries of $G$ have degree at most two,
whereas those of $w$ have degree at most one. Lemma~\ref{lem:rank-one-update}
applies to any positive semidefinite matrix, so for $r=\rank G<n$,
\begin{equation}\label{eq:chart-distance}
 d(L_x,L_y)^2=\frac{e_{r+1}(G+ww\T)}{e_r(G)}.
\end{equation}
For $k<\rho\leqslant r_*$ define
\begin{equation}\label{eq:chart-kernels}
 \widehat F_\rho=e_{\rho+1}(G+ww\T)-e_\rho(G),
 \qquad \widehat R_\rho=e_{\rho+1}(G).
\end{equation}
They are symmetric and satisfy all three assertions of
Proposition~\ref{prop:encoding}, with $G$ in place of $P+Q$.
Their total degrees are at most $2(\rho+1)\leqslant2n$.

Use the degree-$2n$ affine Veronese map on $\R^m$. For a finite
configuration in this chart, the inverse image of the span of its lifted
points is defined by polynomials of degree at most $2n$. The sum of their
squares is a single equation of degree at most $4n$; no equations for a
projection model need to be added. Consequently, Corollary~\ref{lem:uniform-nash}
bounds the number of connected Nash pieces by
\begin{equation}\label{eq:chart-constant}
 C_{k,n}=\Sigma(m,4n)=(32n)^{(4^m-1)/3}.
\end{equation}
The proof of Theorem~\ref{thm:rank-descent} works unchanged on each piece,
using $e_\rho(G)$ as the analytic rank detector and
\eqref{eq:graph-point-distance} at the diagonal.

\begin{proof}[Proof of Theorem~\ref{thm:main}]
Within a fixed chart, let $A^{\mathrm{ch}}_\rho$ be the supremum for
unit-distance families with direction-sum rank at most $\rho$.
The parallel bound and the chart version of rank descent give
\[
 A^{\mathrm{ch}}_k\leqslant n-k+1,\qquad
 A^{\mathrm{ch}}_\rho\leqslant C_{k,n}A^{\mathrm{ch}}_{\rho-1}
 \quad(k<\rho\leqslant r_*).
\]
Hence
\[
 A^{\mathrm{ch}}_{r_*}\leqslant(n-k+1)C_{k,n}^{r_*-k}.
\]
Assign every plane of a finite configuration to one chart containing it.
There are at most $\binom nk$ assigned subfamilies, and
$r_*-k=\min\{k,n-k-1\}=h$. This proves the first inequality in
\eqref{eq:main-bound}. The chart factor is paid once, not at every rank
step.

For the second inequality, the cases $k=0$ and $k=n-1$ follow immediately
from the first bound: it is respectively $n+1$ and $2n$, with $m=n$.
Otherwise $1\leqslant k\leqslant n-2$, so $m\geqslant4$, $n\leqslant m$, and $h\leqslant m$.
Write $B_{k,n}$ for the middle expression in \eqref{eq:main-bound}.
Using $\binom nk\leqslant2^n$, $n-k+1\leqslant2^n$, and
$\log_2(32n)\leqslant5+\log_2m\leqslant2m$, we obtain
\[
 \log_2 B_{k,n}
 \leqslant2m+\frac{h}{3}(4^m-1)\log_2(32n)
 \leqslant2m+\frac23m^2\,4^m.
\]
For integers $m\geqslant4$, $m^2\leqslant2^m$ and
$2m\leqslant\frac13 2^{3m}$. Thus $\log_2B_{k,n}\leqslant2^{3m}$, as required.
\end{proof}

\begin{remark}\label{rem:two-explicit-bounds}
The global-coordinate proof also gives the completely explicit estimate
\begin{equation}\label{eq:global-numeric}
 \mathcal{N}(k,n)\leqslant(n-k+1)
 [16(n+1)]^{\,h(2^{n(n+3)}-1)/3}.
\end{equation}
Both \eqref{eq:global-numeric} and \eqref{eq:main-bound} are valid; the latter
uses the smaller chart dimension. For small or special parameters, the
structural equalities below are much stronger than either expression.
\end{remark}

\section{Embeddings, common factors, and elementary comparisons}\label{sec:basic}

We record the basic operations with their exact parameter ranges.
They will also fix the distinction between identities for all configurations
and identities for a restricted class.

\begin{proposition}[Monotonicity]\label{prop:monotonicity}
For $0\leqslant k<n$ and integers $a,b\geqslant0$,
\begin{equation}\label{eq:monotonicity}
 \mathcal{N}(k,n)\leqslant \mathcal{N}(k+a,n+a+b).
\end{equation}
In particular, $\mathcal{N}(k,n)\leqslant \mathcal{N}(k,n+1)$ and
$\mathcal{N}(k,n)\leqslant \mathcal{N}(k+1,n+1)$.
\end{proposition}
\begin{proof}
Replace each plane by $L\times\R^a$, preserving all distances, and then
embed $\R^{n+a}$ isometrically in $\R^{n+a+b}$.
\end{proof}

\begin{proposition}[Common-direction quotient]\label{prop:quotient}
Suppose that the directions of all members of a family contain a fixed
$h$-dimensional linear space $K$, where $0\leqslant h\leqslant k$. Put
$\ell_L=\proj_{K^\perp}L$. Then
\begin{equation}\label{eq:quotient}
 L=\ell_L+K,\qquad \dim\ell_L=k-h,
 \qquad d(L,L')=d(\ell_L,\ell_{L'}).
\end{equation}
Consequently, the maximum size within this restricted class is exactly
$\mathcal{N}(k-h,n-h)$.
\end{proposition}
\begin{proof}
Each direction decomposes orthogonally as
$U=K\oplus(U\cap K^\perp)$. The component of a base point along $K$ can be
absorbed into the direction. Thus $L=\ell_L+K$, and orthogonality gives
\[
 \norm{(x+u)-(y+v)}^2=\norm{x-y}^2+\norm{u-v}^2
 \quad(x,y\in K^\perp,\ u,v\in K).
\]
Taking infima proves the distance assertion. Conversely, every configuration
in $K^\perp$ can be lifted by adding $K$.
\end{proof}

\begin{corollary}\label{cor:basic-lower}
For all $0\leqslant k<n$,
\begin{equation}\label{eq:parallel-lower}
 \mathcal{N}(k,n)\geqslant n-k+1.
\end{equation}
For $1\leqslant k<n$,
\begin{equation}\label{eq:line-lower}
 \mathcal{N}(k,n)\geqslant \mathcal{N}(1,n-k+1).
\end{equation}
\end{corollary}
\begin{proof}
Lift a regular simplex of side length one by a common $k$-dimensional
factor for the first assertion, and lift a line configuration by a common
$(k-1)$-dimensional factor for the second.
\end{proof}

\section{Direction dichotomy, sharp top bounds, and exact reductions}
\label{sec:reductions}

\subsection{A direction-space dichotomy}
\begin{lemma}[Star--top dichotomy]\label{lem:star-top}
Let $\{U_i\}$ be a nonempty family of $k$-dimensional linear subspaces,
where $1\leqslant k<n$, satisfying
\begin{equation}\label{eq:codim-one-intersections}
 \dim(U_i\cap U_j)\geqslant k-1\qquad\text{for all }i,j.
\end{equation}
Then at least one of the following holds:
\begin{enumerate}
\item all $U_i$ contain a common $(k-1)$-dimensional space;
\item all $U_i$ lie in a common $(k+1)$-dimensional space.
\end{enumerate}
\end{lemma}
\begin{proof}
If $k=1$, the first alternative holds with the zero space. If all directions
coincide, it also holds. Otherwise choose distinct $U_1,U_2$ and put
$K=U_1\cap U_2$, $W=U_1+U_2$; these have dimensions $k-1$ and $k+1$.
If every $U_i$ contains $K$, we are done. Choose $U_3$ not containing $K$.
The spaces $U_3\cap U_1$ and $U_3\cap U_2$ are distinct hyperplanes of
$U_3$. They span $U_3$, so $U_3\subset W$ and
$\dim(U_3\cap K)=k-2$.

The same reasoning puts any $U_i$ not containing $K$ inside $W$.
If some $U_i$ contains $K$ but is not contained in $W$, then
$U_i\cap W=K$, whence
$\dim(U_i\cap U_3)=k-2$, contradicting
\eqref{eq:codim-one-intersections}. Therefore every $U_i\subset W$.
\end{proof}

This is a linear-algebraic statement, independent of the Euclidean metric;
compare \cite[Proposition~17 and Corollary~18]{GR}. In an equidistant
configuration, the condition $\dim(U_i+U_j)\leqslant k+1$ is exactly
\eqref{eq:codim-one-intersections}. Repeated directions cause no problem,
since the lemma can first be applied to the set of distinct directions.

\subsection{Distances in the top case}
Fix a $(k+1)$-dimensional space $W\subset\R^n$, and put
\begin{equation}\label{eq:top-split}
 c=n-k,\qquad E=W^\perp,\qquad\dim E=c-1.
\end{equation}
If every direction lies in $W$, each affine plane has the form
$L_i=z_i+H_i$, where $z_i\in E$ and $H_i$ is an affine hyperplane of $W$.
Partition the family into parallel-direction classes, indexed by $a$.
For a class $a$, choose a common unit normal $u_a\in W$ and write
\[
 H_i=\{w\in W|u_a\cdot w=t_i\}\qquad(i\text{ in class }a).
\]
Distinct classes have nonparallel hyperplanes in $W$, hence intersect
in that factor. Thus
\begin{equation}\label{eq:top-distances}
 d(L_i,L_j)^2=
 \begin{cases}
  \norm{z_i-z_j}^2+(t_i-t_j)^2,&i,j\text{ in the same class},\\
  \norm{z_i-z_j}^2,&i,j\text{ in different classes}.
 \end{cases}
\end{equation}

\begin{theorem}[Sharp top bound]\label{thm:top}
For $1\leqslant k<n$, let $\cT(k,n)$ be the maximum size of an equidistant family
whose directions all lie in a fixed $(k+1)$-dimensional space. Then
\begin{equation}\label{eq:top-bound}
 \cT(k,n)=2(n-k).
\end{equation}
\end{theorem}

\begin{proof}
We use \eqref{eq:top-split}--\eqref{eq:top-distances}. If there is only one
direction class, the parallel bound gives at most $c+1\leqslant2c$ members.
Assume that there are at least two classes. For each class $a$, let $m_a$
be its number of members, let $Z_a$ be the set of its distinct $E$-coordinates,
and let $E_a$ be the direction space of $\aff Z_a$, of dimension $d_a$.
The pairs $(z_i,t_i)$ in this class are equidistant in an affine space
of dimension at most $d_a+1$, so
\begin{equation}\label{eq:class-size}
 m_a\leqslant d_a+2.
\end{equation}

All cross-distances between $Z_a$ and $Z_b$ are one when $a\ne b$.
Subtracting four squared-distance equations gives
$E_a\perp E_b$. Fix a point $z$ in another class. Its orthogonal projection
onto $\aff Z_a$ is a center for a sphere containing $Z_a$.
A sphere center within $\aff Z_a$ is unique: subtracting two equal-distance
equations determines its scalar products with the spanning differences
of $Z_a$. Denote this center by $o_a$ and its radius by $r_a$.
Every point of every other class projects to the same $o_a$.
By affine combinations, the same holds for the other centers. It follows
that
\begin{equation}\label{eq:orth-centers}
 E_a\perp E_b\ (a\ne b),\qquad
 o_a-o_b\perp E_j\quad\text{for all }a,b,j,
\end{equation}
and orthogonal expansion of a cross-distance gives
\begin{equation}\label{eq:center-distance}
 \norm{o_a-o_b}^2+r_a^2+r_b^2=1\qquad(a\ne b).
\end{equation}

Let $J=\{a|r_a^2<1/2\}$ and $q=|J|$. If $d_a=0$, then $r_a=0$.
If $d_a=1$, the set $Z_a$ lies on a sphere in a line and spans that line,
so it consists of the two points at distance $2r_a$.
Within one class \eqref{eq:top-distances} gives $\norm{z_i-z_j}\leqslant1$,
and hence $r_a\leqslant1/2$. Consequently $a\notin J$ implies $d_a\geqslant2$.
From \eqref{eq:class-size} we therefore obtain
\begin{equation}\label{eq:size-by-radii}
 \sum_a m_a\leqslant2\sum_a d_a+2q.
\end{equation}
Indeed, for $a\notin J$ use $d_a+2\leqslant2d_a$, and for $a\in J$ use
$d_a+2\leqslant2d_a+2$.

If $q\geqslant1$, the centers $o_a$ with $a\in J$ are affinely independent.
For coefficients with $\sum\limits_{a\in J}\lambda_a=0$,
\eqref{eq:center-distance} yields
\begin{equation}\label{eq:center-affine-independence}
 \left\|\sum_{a\in J}\lambda_a o_a\right\|^2
 =-\frac12\sum_{a,b\in J}\lambda_a\lambda_b\norm{o_a-o_b}^2
 =\sum_{a\in J}\left(\frac12-r_a^2\right)\lambda_a^2.
\end{equation}
The last expression is strictly positive for nonzero $\lambda$.
By \eqref{eq:orth-centers}, their affine direction space is orthogonal
to the mutually orthogonal spaces $E_a$. Therefore
\[
 c-1=\dim E\geqslant\sum_a d_a+q-1.
\]
Together with \eqref{eq:size-by-radii}, this gives $\sum\limits_a m_a\leqslant2c$.
If $q=0$, the stronger estimate $\sum\limits_a m_a\leqslant2\sum\limits_a d_a\leqslant2(c-1)$
holds. This proves the upper bound.

For the lower bound, choose the $c$ vertices $z_1,\ldots,z_c$ of a
regular simplex of side length one in $E$. For each $a$, choose a distinct
hyperplane direction in $W$ and two parallel affine hyperplanes
$H_a^+,H_a^-$ of that direction at distance one. Arbitrarily many distinct
hyperplane directions exist because $\dim W=k+1\geqslant2$.
The $2c$ planes $z_a+H_a^\pm$ have all mutual distances one by
\eqref{eq:top-distances}. For $c=1$, take a single simplex vertex and
one parallel pair. This proves equality.
\end{proof}

\subsection{Exact reduction at direction-sum rank \texorpdfstring{$k+1$}{k+1}}
\begin{proof}[Proof of \eqref{eq:rank-one-reduction} in Theorem~\ref{thm:reductions}]
Lifting a line configuration in $\R^{n-k+1}$ by a common
$(k-1)$-dimensional factor preserves distance and gives direction-sum
rank at most $k+1$. This proves the lower bound.
For the reverse inequality, use Lemma~\ref{lem:star-top}.
In the star case, Proposition~\ref{prop:quotient} reduces the family to
lines in $\R^{n-k+1}$.

In the top case, use \eqref{eq:top-distances}. Choose, for each direction
class $a$, a unit vector $v_a\in\R^2$, with no two $v_a$ parallel.
Replace $L_i=z_i+H_i$ in class $a$ by the line
\begin{equation}\label{eq:top-compression}
 \ell_i=\{x\in\R^2|v_a\cdot x=t_i\}\times\{z_i\}
 \subset\R^2\oplus E.
\end{equation}
Within a class, the distance in the planar factor is $|t_i-t_j|$;
between classes, the planar lines intersect. Thus \eqref{eq:top-distances}
shows that every pairwise distance is preserved. The ambient dimension
is $2+(n-k-1)=n-k+1$, giving the claimed upper bound.
For $k=1$, the conclusion is also immediate from the definition.
\end{proof}

\begin{corollary}[Large rank-one-excess families split]\label{cor:large-star}
If $1\leqslant k<n$, $\dim(U_i+U_j)\leqslant k+1$, and an equidistant family has more than
$2(n-k)$ members, then its directions contain a common $(k-1)$-dimensional
subspace.
\end{corollary}
\begin{proof}
The top alternative in Lemma~\ref{lem:star-top} is excluded by
Theorem~\ref{thm:top}.
\end{proof}

\begin{corollary}\label{cor:two-c}
For $1\leqslant k<n$,
\begin{equation}\label{eq:two-c}
 \mathcal{N}(k,n)\geqslant2(n-k).
\end{equation}
If $n-k\geqslant2$, then
$\mathcal{N}(k,n)\geqslant\max\{\mathcal{N}(1,3),2(n-k)\}$.
\end{corollary}
\begin{proof}
Use the construction in Theorem~\ref{thm:top}; for the other term lift
three-dimensional lines by $\R^{k-1}$ and then embed in $\R^n$.
\end{proof}

\subsection{Codimension two: an actual common factor}
\begin{proof}[Proof of \eqref{eq:codim-splitting}--\eqref{eq:codim-equality}]
Write $k=n-2\geqslant1$. Positive separation implies
$\dim(U_i+U_j)\leqslant n-1=k+1$, so Lemma~\ref{lem:star-top} applies.
The star alternative immediately supplies a common $(n-3)$-space.
In the top alternative, all directions are hyperplanes of a common
$(n-1)$-space $W$. In \eqref{eq:top-split}, $E$ is now a line.
By \eqref{eq:top-distances}, representatives of distinct direction classes
have pairwise distance one in $E$. There are therefore at most two
classes. If there is one direction, choose any $(n-3)$-subspace of it.
If there are two, their intersection in $W$ has dimension $n-3$.
In either case all directions contain a common $(n-3)$-space $K$.
Proposition~\ref{prop:quotient} gives \eqref{eq:codim-splitting} and the
upper bound in \eqref{eq:codim-equality}. The reverse bound follows by
$\ell\mapsto\ell\times\R^{n-3}$.
\end{proof}

This proves more than equality of extremal numbers: every codimension-two
configuration has the stated product form, regardless of its cardinality.
For $n=3$, the common factor is the zero space. The reduction itself uses
neither H\"ofer's bound nine nor a computer-assisted bound seven.

\begin{corollary}[Specialization using the three-dimensional theorem]
\label{cor:codim-seven}
Using the computer-assisted upper bound of
\cite[Theorem~1.1]{MRS} and the seven-line construction of \cite{BLR},
\begin{equation}\label{eq:all-codim-seven}
 \mathcal{N}(n-2,n)=7\qquad(n\geqslant3).
\end{equation}
\end{corollary}
\begin{proof}
The two cited results give $\mathcal{N}(1,3)=7$. Apply
\eqref{eq:codim-equality}.
\end{proof}

\begin{remark}[External input]\label{rem:external-seven}
Corollary~\ref{cor:codim-seven} explicitly imports the three-dimensional
nonexistence theorem from the public preprint \cite{MRS}; no new
computer-free proof of that theorem is asserted here.
Without using that input, \cite{BLR,Hoefer} and the same reduction give
$7\leqslant \mathcal{N}(n-2,n)\leqslant9$. All parameter identities in this section hold
independently of the numerical specialization.
\end{remark}
\section{Transverse products}\label{sec:products}

For nonempty affine subspaces in two orthogonal factors, independence of
the two minimizations gives
\begin{equation}\label{eq:distance-product}
 d(A\times B,A'\times B')^2=d(A,A')^2+d(B,B')^2.
\end{equation}
The ordinary product of two unit-distance configurations therefore has
some distances $\sqrt2$. The next lemma makes the second contribution
vanish whenever the first is nonzero.

\begin{lemma}[Simultaneous transverse rotations]\label{lem:rotations}
Let $V_1,\ldots,V_q$ be $k$-dimensional subspaces of $\R^n$, with
$0<k<n$ and $2k\geqslant n$. For every positive integer $p$ there exist
$g_1,\ldots,g_p\in SO(n)$ such that
\begin{equation}\label{eq:transversality}
 g_iV_j+g_{i'}V_{j'}=\R^n
 \qquad(i\ne i',\ 1\leqslant j,j'\leqslant q).
\end{equation}
\end{lemma}
\begin{proof}
For fixed $k$-spaces $U,V$, the set
$\left\{g\in SO(n)\big|U+gV=\R^n\right\}$ is open. It is nonempty because
$2k\geqslant n$: a $k$-space can be chosen to contain $U^\perp$, and $SO(n)$
acts transitively on the unoriented $k$-planes. In matrix coordinates,
the failure of the condition is the simultaneous vanishing of the
$n\times n$ minors of concatenated bases for $U$ and $gV$.
The sum of their squares is a real-analytic function on connected $SO(n)$
that is not identically zero. Its zero set has empty interior, so the
set of good rotations is dense. Here $n\geqslant2$, as follows from the
hypotheses.

Choose the rotations successively. At each stage there are only finitely
many pairs of previous spaces and current spaces to make transverse.
The intersection of their open dense sets of good rotations is nonempty.
This proves \eqref{eq:transversality}.
\end{proof}

\begin{proof}[Proof of Theorem~\ref{thm:product}]
Let $A_1,\ldots,A_p$ be an equidistant $k_1$-plane configuration in
$\R^{n_1}$, and let $B_1,\ldots,B_q$ be an equidistant $k_2$-plane
configuration in $\R^{n_2}$, with directions $V_j$.
Assume $2k_2\geqslant n_2$ and choose the rotations of
Lemma~\ref{lem:rotations}. Define
\begin{equation}\label{eq:product-construction}
 L_{i,j}=A_i\times g_iB_j
 \subset\R^{n_1}\oplus\R^{n_2}.
\end{equation}
If $i=i'$ and $j\ne j'$, \eqref{eq:distance-product} gives squared distance
$0+1$. If $i\ne i'$, the affine spaces $g_iB_j$ and $g_{i'}B_{j'}$
intersect: their directions sum to the whole second factor, so every
translation difference can be absorbed into their directions.
Their separation is therefore zero, and the squared distance in
\eqref{eq:distance-product} is $1+0$. Thus all $pq$ members of
\eqref{eq:product-construction} are distinct and at mutual distance one.
Finiteness has already been proved, so the integer-valued suprema
$\mathcal{N}(k_i,n_i)$ are attained. Take $p=\mathcal{N}(k_1,n_1)$ and $q=\mathcal{N}(k_2,n_2)$.
If the condition holds for the first factor instead, interchange the two
factors.
\end{proof}

\begin{corollary}[Padded product]\label{cor:padded-product}
For arbitrary $0\leqslant k_i<n_i$, define
\begin{equation}\label{eq:padding}
 t=\min\{(n_1-2k_1)_+,(n_2-2k_2)_+\},\qquad x_+=\max\{x,0\}.
\end{equation}
Then
\begin{equation}\label{eq:padded-product}
 \mathcal{N}(k_1+k_2+t,n_1+n_2+t)
 \geqslant \mathcal{N}(k_1,n_1)\mathcal{N}(k_2,n_2).
\end{equation}
\end{corollary}
\begin{proof}
Select a factor attaining the minimum in \eqref{eq:padding} and add
$\R^t$ to every member of that configuration. Its parameters become
$(k_i+t,n_i+t)$, with $2(k_i+t)\geqslant n_i+t$.
Proposition~\ref{prop:monotonicity} preserves its cardinality, and
Theorem~\ref{thm:product} applies to the enlarged factor and the other one.
\end{proof}

\begin{remark}
The transversality condition concerns different rotated copies of one
configuration. It does not require two members within the same copy to
have transverse directions; that would force their separation to be
zero. The independent rotations in \eqref{eq:product-construction} are
therefore essential.
\end{remark}

\section{Consequences for the parameter array}\label{sec:consequences}

\subsection{A supermultiplicative sequence}
Recall $S(c)=\mathcal{N}(c-1,2c-1)$ for $c\geqslant1$ and $a=\mathcal{N}(1,3)$.

\begin{proposition}\label{prop:supermultiplicative}
For positive integers $c,d$,
\begin{equation}\label{eq:supermultiplicative}
 S(c+d)\geqslant S(c)S(d).
\end{equation}
Moreover,
\begin{equation}\label{eq:S-lower}
 S(c)\geqslant2^{c\bmod2}a^{\lfloor c/2\rfloor}.
\end{equation}
\end{proposition}
\begin{proof}
Increase the dimension of every member of a configuration counted by
$S(d)$ by one common direction. This gives
$\mathcal{N}(d,2d)\geqslant \mathcal{N}(d-1,2d-1)=S(d)$. The second factor now satisfies the
transversality inequality with equality, and
Theorem~\ref{thm:product} gives
\[
 \mathcal{N}(c+d-1,2c+2d-1)
 \geqslant \mathcal{N}(c-1,2c-1)\mathcal{N}(d,2d)\geqslant S(c)S(d).
\]
Finally $S(1)=\mathcal{N}(0,1)=2$ and $S(2)=a$. Iterating with blocks of size two,
and one additional block of size one when needed, proves
\eqref{eq:S-lower}.
\end{proof}

\begin{corollary}[Lower bounds at fixed codimension]\label{cor:fixed-codim-lower}
If $c\geqslant1$ and $k\geqslant c-1$, then
\begin{equation}\label{eq:fixed-codim-lower}
 \mathcal{N}(k,k+c)\geqslant2^{c\bmod2}a^{\lfloor c/2\rfloor}
 \geqslant2^{c\bmod2}7^{\lfloor c/2\rfloor}.
\end{equation}
In particular,
\begin{equation}\label{eq:growth}
 \liminf_{c\to\infty}\mathcal{N}(c-1,2c-1)^{1/c}\geqslant\sqrt7.
\end{equation}
\end{corollary}
\begin{proof}
Add a common $(k-c+1)$-dimensional factor to a configuration counted by
$S(c)$, and apply \eqref{eq:S-lower} and the lower bound $a\geqslant7$ from
\cite{BLR}. Taking $c$th roots gives \eqref{eq:growth}.
\end{proof}

The numerical examples in Table~\ref{tab:examples} use only the certified
seven-line construction, not the upper bound in \cite{MRS}.
\begin{table}[htbp]
\centering
\caption{Lower bounds supplied by transverse products.}
\label{tab:examples}
\begin{tabular}{@{}ccc@{}}
\toprule
$(k,n)$ & Symbolic lower bound & Numerical lower bound\\
\midrule
$(2,5)$ & $2a$ & $14$\\
$(3,7)$ & $a^2$ & $49$\\
$(4,9)$ & $2a^2$ & $98$\\
$(5,11)$ & $a^3$ & $343$\\
$(7,15)$ & $a^4$ & $2401$\\
\bottomrule
\end{tabular}
\end{table}

\subsection{Simplex amplification}
\begin{proposition}\label{prop:simplex-amplification}
If $1\leqslant k<q\leqslant2k$ and $n\geqslant q$, then
\begin{equation}\label{eq:simplex-amplification}
 \mathcal{N}(k,n)\geqslant(n-q+1)\mathcal{N}(k,q).
\end{equation}
Consequently,
\begin{align}
 \mathcal{N}(k,n)&\geqslant(n-k-1)a
 &&(k\geqslant2,\ n\geqslant k+2),\label{eq:amplify-a}\\
 \mathcal{N}(k,n)&\geqslant(n-2k+1)\,2^{k\bmod2}a^{\lfloor k/2\rfloor}
 &&(k\geqslant1,\ n\geqslant2k).\label{eq:amplify-balanced}
\end{align}
\end{proposition}
\begin{proof}
If $n=q$, \eqref{eq:simplex-amplification} is tautological.
Otherwise take a regular simplex of $n-q+1$ points in $\R^{n-q}$
as the first factor of Theorem~\ref{thm:product} and a configuration
counted by $\mathcal{N}(k,q)$ as the second factor. The hypothesis $q\leqslant2k$
ensures transversality. For \eqref{eq:amplify-a}, take $q=k+2$ and use
codimension-two reduction. For \eqref{eq:amplify-balanced}, take $q=2k$,
observe that $\mathcal{N}(k,2k)\geqslant \mathcal{N}(k-1,2k-1)=S(k)$ by a common-factor lift,
and apply \eqref{eq:S-lower}.
\end{proof}

For $k=1$, \eqref{eq:amplify-balanced} recovers the lower bound
$2n-2$ of \cite[Section~3]{DKL}. The other cases amplify
higher-dimensional configurations in the same way. Different lower
bounds in this paper may of course be combined by taking their maximum.

\subsection{Stopping the upper-bound recursion earlier}
The exact reduction at rank $k+1$ can be substituted into the upper-bound
argument, linking the structural theorem with the quantitative proof.

\begin{corollary}\label{cor:early-stop}
Suppose $1\leqslant k\leqslant n-2$, and retain $m,h$ from \eqref{eq:m-h}.
Then
\begin{equation}\label{eq:early-stop}
 \mathcal{N}(k,n)\leqslant\binom nk\,\mathcal{N}(1,n-k+1)
 (32n)^{\,(h-1)(4^m-1)/3}.
\end{equation}
Writing $c=n-k$, the external line bound \eqref{eq:SZ} gives
\begin{equation}\label{eq:early-stop-explicit}
 \mathcal{N}(k,n)\leqslant\binom nk\,4(c+1)7^{2c-1}
 (32n)^{\,(h-1)(4^m-1)/3}.
\end{equation}
\end{corollary}
\begin{proof}
In a fixed chart, start the recursion with
$A^{\mathrm{ch}}_{k+1}\leqslant \mathcal{N}_{\leqslant k+1}(k,n)=\mathcal{N}(1,n-k+1)$.
There are $r_*-(k+1)=h-1$ remaining descents.
Multiply by the chart count once, as in Theorem~\ref{thm:main}.
Since $c\geqslant2$, \eqref{eq:SZ} is applicable with ambient dimension $c+1$.
\end{proof}

Within the stated parameter range, the external line bound in
\eqref{eq:early-stop-explicit} makes this estimate strictly smaller than
the chart bound in \eqref{eq:main-bound}. When $c=2$, the exact equality
\eqref{eq:codim-equality} is substantially stronger.

\subsection{Fixed-codimension monotonicity and its limitations}
For $c\geqslant1$ put $F_c(k)=\mathcal{N}(k,k+c)$, $k\geqslant0$.
Proposition~\ref{prop:monotonicity} gives
\begin{equation}\label{eq:fixed-codim-monotone}
 F_c(k+1)\geqslant F_c(k).
\end{equation}
We have $F_1(k)=2$ for all $k\geqslant0$, while
$F_2(k)=a$ for $k\geqslant1$ and $F_2(0)=3$.
A universal common-factor theorem does not extend from codimension two
to codimension three, as the following example shows.

\begin{example}[No forced common direction in codimension three]
\label{ex:codim-three}
In $\R^5$ let
\begin{equation}\label{eq:codim-three-example}
 \begin{split}
 L_1&=\Span(e_1,e_2),\\
 L_2&=e_5+\Span(e_3,e_4),\\
 L_3&=e_2+e_4+\Span(e_1+e_3,e_5).
 \end{split}
\end{equation}
The directions $U_1,U_2,U_3$ have pairwise zero intersection.
Their pairwise sums have orthogonal complements
\[
 (U_1+U_2)^\perp=\R e_5,\quad
 (U_1+U_3)^\perp=\R e_4,\quad
 (U_2+U_3)^\perp=\R e_2.
\]
The corresponding translation differences have projection of length one
onto these lines, so all three separation distances are one.
In particular, the configuration has no common nonzero direction.
\end{example}

This example rules out a direct extension of the common-factor statement,
not a possible stabilization of extremal numbers. We do not assert
$\mathcal{N}(k,k+c)=\mathcal{N}(c-1,2c-1)$ for $c\geqslant3$ and $k\geqslant c-1$; only the lower-bound
direction follows from adding a common factor.

\section{Fixed-codimension bounds and stabilization}\label{sec:stability}

We study $F_c(k)=\mathcal{N}(k,k+c)$ for fixed codimension $c$.
Monotonicity shows that this sequence is nondecreasing. We prove that
it is bounded, and hence eventually constant, by projecting to a
uniformly bounded span of normal spaces.

\subsection{A distance-preserving projection lemma}

Let $L_i=p_i+U_i$ be affine subspaces of $\R^n$, and put
$V_i=U_i^\perp$. The distance formula gives
\begin{equation}\label{eq:normal-distance}
d(L_i,L_j)=
\left\|\operatorname{proj}_{V_i\cap V_j}(p_i-p_j)\right\|.
\end{equation}

\begin{lemma}[Normal-space projection]\label{lem:normal-projection}
Let $W\subset\R^n$ be linear and let $M_i=\operatorname{proj}_W(L_i)$.
The normal space of $M_i$ inside $W$ is $V_i\cap W$, and
\begin{equation}\label{eq:projected-distance}
d_W(M_i,M_j)=
\left\|\operatorname{proj}_{W\cap V_i\cap V_j}(p_i-p_j)\right\|.
\end{equation}
In particular, if $V_i\cap V_j\subset W$, then projection to $W$
preserves the distance of this pair. A sufficient condition is
$\pi_W(V_i)\cap\pi_W(V_j)=\{0\}$, where
$\pi_W:\R^n\to\R^n/W$ is the quotient map.
\end{lemma}

\begin{proof}
The direction of $M_i$ is $\operatorname{proj}_W(U_i)$. For $w\in W$,
$w\perp\operatorname{proj}_W(U_i)$ if and only if $w\perp U_i$.
Thus its normal space in $W$ is $W\cap V_i$. The distance formula
inside $W$ gives \eqref{eq:projected-distance}. If the quotient
images intersect trivially, every vector in $V_i\cap V_j$ maps to
zero, and hence belongs to $W$.
\end{proof}

\subsection{A uniform normal-span bound in each codimension}

For $c\geqslant1$, put
\[
b_c=\sup_{k\geqslant0}\mathcal{N}(k,k+c),
\]
initially allowing $b_c=+\infty$.

\begin{theorem}[Uniform normal-span bound]\label{thm:normal-span}
For every $c\geqslant1$, there is an integer $D_c\geqslant c$ such
that every finite unit-distance family of affine subspaces of
codimension $c$ satisfies
\[
\dim\left(\sum_i U_i^\perp\right)\leqslant D_c.
\]
One may take
\begin{equation}\label{eq:Dc-recursion}
D_1=1,\qquad
D_c=c+\sum_{t=1}^{c-1}t(c-t)b_t\quad(c\geqslant2),
\end{equation}
where the lower-codimension numbers $b_t$ are finite by induction.
\end{theorem}

\begin{proof}
For $c=1$, positive separation forces all affine hyperplanes to be
parallel. Thus $D_1=1$ and $b_1=2$.

Assume $c\geqslant2$ and that $b_t<\infty$ for $1\leqslant t<c$.
Consider a finite unit-distance family $L_i=p_i+U_i$ in $\R^n$,
and write $V_i=U_i^\perp$, so $\dim V_i=c$. Start with $W_1=V_1$.
We construct nested spaces $W_1\subset\cdots\subset W_c$ such that
\begin{equation}\label{eq:intersection-growth}
\dim(V_i\cap W_r)\geqslant r\qquad\text{for every }i.
\end{equation}
For $r=1$, this follows from positive separation and
\eqref{eq:normal-distance}.

Suppose $W_r$ has been constructed, where $r<c$. Let
$\pi_r:\R^n\to\R^n/W_r$ be the quotient map. Among indices with
$V_i\not\subset W_r$, choose a maximal subfamily
$V_{i_1},\ldots,V_{i_s}$ whose nonzero quotient images are pairwise
disjoint (their pairwise intersections are zero). Put
$t_\nu=\dim(V_{i_\nu}\cap W_r)$. Then
$r\leqslant t_\nu\leqslant c-1$.

For each fixed $t$, the projections to $W_r$ of the selected
affine subspaces with $t_\nu=t$ have codimension $t$. Their quotient
normal spaces are pairwise disjoint, so Lemma~\ref{lem:normal-projection}
shows that the projected subspaces remain pairwise at distance one.
They are distinct, and therefore
\[
\#\{\nu|t_\nu=t\}\leqslant
\mathcal{N}(\dim W_r-t,\dim W_r)\leqslant b_t.
\]
Set
\[
W_{r+1}=W_r+\sum_{\nu=1}^sV_{i_\nu}.
\]
Each selected space adds at most $c-t_\nu$ dimensions, whence
\[
\dim W_{r+1}-\dim W_r
\leqslant\sum_{t=r}^{c-1}(c-t)b_t.
\]

We verify \eqref{eq:intersection-growth} at level $r+1$.
It holds for selected spaces because they are contained in
$W_{r+1}$. If $V_i\not\subset W_r$ is not selected, maximality
gives a selected $V_{i_\nu}$ such that
$\pi_r(V_i)\cap\pi_r(V_{i_\nu})\ne\{0\}$. Hence there are
$v\in V_i\setminus W_r$ and $w\in V_{i_\nu}$ with $v-w\in W_r$.
Thus
\[
v\in V_i\cap(W_r+V_{i_\nu})\subset V_i\cap W_{r+1},
\qquad v\notin W_r.
\]
Together with the old $r$-dimensional intersection, this proves
the desired increase. Spaces already contained in $W_r$ require
no argument.

After at most $c-1$ steps, every $V_i$ is contained in $W_c$.
Summing the dimension increments gives
\[
\begin{aligned}
\dim W_c
&\leqslant c+\sum_{r=1}^{c-1}\sum_{t=r}^{c-1}(c-t)b_t\\
&=c+\sum_{t=1}^{c-1}t(c-t)b_t=D_c.
\end{aligned}
\]
\end{proof}

\begin{theorem}[Fixed-codimension finiteness and stabilization]
\label{thm:fixed-codim-stability}
For every $c\geqslant1$, the number
$b_c=\sup\limits_{k\geqslant0}\mathcal{N}(k,k+c)$ is finite. With $D_c$
as in Theorem~\ref{thm:normal-span}, one has
\begin{equation}\label{eq:fixed-codim-uniform}
\mathcal{N}(k,k+c)\leqslant\mathcal{N}(D_c-c,D_c)
\qquad(k\geqslant0),
\end{equation}
and
\begin{equation}\label{eq:fixed-codim-stable}
\mathcal{N}(k,k+c)=\mathcal{N}(D_c-c,D_c)
\qquad(k\geqslant D_c-c).
\end{equation}
In particular, $b_c=\mathcal{N}(D_c-c,D_c)$.
\end{theorem}

\begin{proof}
We argue by induction on $c$. Theorem~\ref{thm:normal-span} gives
$D_c<\infty$ once $b_1,\ldots,b_{c-1}$ are finite. For a finite
unit-distance family of affine $(n-c)$-planes in $\R^n$, let
\[
W=\sum_iU_i^\perp,\qquad d=\dim W\leqslant D_c.
\]
Every direction $U_i$ contains $W^\perp$. Quotienting by this common
direction gives a unit-distance family of affine $(d-c)$-planes in
$W\cong\R^d$. Adding a common $(D_c-d)$-dimensional factor embeds
it into $\R^{D_c}$ as affine $(D_c-c)$-planes, preserving distances.
Its cardinality is therefore at most $\mathcal{N}(D_c-c,D_c)$.
This proves \eqref{eq:fixed-codim-uniform}. The right-hand side is
finite by Theorem~\ref{thm:main}, so $b_c$ is finite.

If $k\geqslant D_c-c$, monotonicity gives
$\mathcal{N}(D_c-c,D_c)\leqslant\mathcal{N}(k,k+c)$.
Together with \eqref{eq:fixed-codim-uniform}, this proves
\eqref{eq:fixed-codim-stable}. The induction starts with $c=1$,
where $D_1=1$ and $\mathcal{N}(0,1)=2$.
\end{proof}

\begin{remark}
The recursion \eqref{eq:Dc-recursion} is an explicit dimension
recursion in terms of the already finite lower-codimension extremal
numbers. It is not, by itself, a closed numerical formula for $D_c$.
\end{remark}

\subsection{A sharper normal-span bound in codimension three}

Put $a=\mathcal{N}(1,3)$. We use the codimension-two splitting theorem
and the sharp top bound to improve the general recursion when $c=3$.

\begin{theorem}[Codimension-three normal-span bound]
\label{thm:codim-three-normal}
Every finite unit-distance family of affine subspaces of codimension
three satisfies
\begin{equation}\label{eq:codim-three-normal}
\dim\left(\sum_iU_i^\perp\right)\leqslant a+3.
\end{equation}
\end{theorem}

\begin{proof}
Write $V_i=U_i^\perp$, so $\dim V_i=3$. If
$\dim(V_i\cap V_j)\geqslant2$ for every pair, apply the
star--top dichotomy to $\{V_i\}$. In the top case,
$\dim\sum_iV_i\leqslant4$. In the star case, all $V_i$ contain
a common two-dimensional space $K$. Distinct three-dimensional
spaces then intersect exactly in $K$. Projecting the affine
subspaces to $K$ gives points at mutual distance one, so there
are at most three distinct direction classes. Hence
$\dim\sum_iV_i\leqslant2+3=5$.

We may therefore assume that some pair, relabelled $V_0,V_1$,
satisfies $V_0\cap V_1=E$ with $\dim E=1$. Put $W=V_0+V_1$,
so $\dim W=5$. We claim that
\begin{equation}\label{eq:codim3-intersection-W}
\dim(V_i\cap W)\geqslant2\qquad\text{for every }i.
\end{equation}
Otherwise the positive intersections with $V_0$ and $V_1$
would both equal $E$. If $e$ is a unit vector spanning $E$ and
$p_i$ is a point on the corresponding affine subspace, then
\eqref{eq:normal-distance} gives
\[
|p_0\cdot e-p_1\cdot e|
=|p_0\cdot e-p_i\cdot e|
=|p_1\cdot e-p_i\cdot e|=1,
\]
which is impossible for three real numbers.

For each $V_i\not\subset W$, its image in $\R^n/W$ is a line,
by \eqref{eq:codim3-intersection-W}. Choose one representative
from each distinct quotient line, say $V_{i_1},\ldots,V_{i_q}$,
and put $A_j=V_{i_j}\cap W$. Then $\dim A_j=2$, and the quotient
lines are pairwise disjoint. Lemma~\ref{lem:normal-projection}
shows that the projections to $W$ of the corresponding affine
subspaces form a unit-distance configuration of affine
three-planes in $W\cong\R^5$, with normal spaces $A_j$.

The codimension-two splitting theorem, applied in $W$, gives a
three-dimensional space $Z\subset W$ containing every $A_j$.
If $q\leqslant4$, then
\[
\dim\sum_iV_i\leqslant5+q\leqslant9\leqslant a+3,
\]
since $a\geqslant7$. We may assume $q\geqslant5$.

The spaces $A_j$ have no common nonzero vector. Indeed, otherwise
all directions $A_j^\perp$ of the projected three-planes would
lie in a fixed four-dimensional subspace of $W$, and the sharp
top bound would imply $q\leqslant2(5-3)=4$.

Set $R_0=V_0\cap Z$ and $R_1=V_1\cap Z$. Each $A_j$ intersects
both $V_0$ and $V_1$ nontrivially, so it intersects both $R_0$
and $R_1$ nontrivially. If $\dim R_0=1$, every $A_j$ would contain
that line, a contradiction. Thus $\dim R_0\geqslant2$, and
similarly $\dim R_1\geqslant2$. Neither dimension is three:
if $R_0=Z$, then $R_1=V_1\cap Z=V_1\cap V_0=E$ would have
dimension one. Hence $\dim R_0=\dim R_1=2$. Since both lie in
the three-dimensional space $Z$, their intersection is nonzero;
it is contained in $V_0\cap V_1=E$. Therefore
\begin{equation}\label{eq:R0R1-intersection}
R_0\cap R_1=E.
\end{equation}

Let $p_b$ be a point on the original affine subspace with normal
space $V_b$, for $b=0,1$, and define affine three-planes in $W$ by
\[
\widetilde M_b=\operatorname{proj}_W(p_b)+R_b^{\perp_W}.
\]
For every $j$,
\[
R_b\cap A_j=(V_b\cap Z)\cap A_j
=V_b\cap A_j=V_b\cap V_{i_j},
\]
where $V_b\subset W$ and $A_j\subset Z\subset W$. Thus the normal
intersections determining the distances from $\widetilde M_b$
to the projected members are exactly the original normal
intersections. Formula \eqref{eq:normal-distance} shows that
\[
d_W(\widetilde M_b,M_{i_j})=1.
\]
Likewise, \eqref{eq:R0R1-intersection} gives
\[
R_0\cap R_1=V_0\cap V_1,
\qquad
d_W(\widetilde M_0,\widetilde M_1)=1.
\]
Consequently, the $q+2$ affine three-planes
\[
M_{i_1},\ldots,M_{i_q},\widetilde M_0,\widetilde M_1
\]
form a unit-distance configuration in $\R^5$. The codimension-two
identity gives $\mathcal{N}(3,5)=\mathcal{N}(1,3)=a$, so
$q+2\leqslant a$.

Every normal space not contained in $W$ has quotient image equal
to one of the selected quotient lines. It is therefore contained
in $W+V_{i_j}$ for some $j$. Hence
\[
\dim\sum_iV_i\leqslant5+q\leqslant a+3,
\]
as required.
\end{proof}

\begin{corollary}[Codimension-three stabilization]\label{cor:codim-three-stable}
For every $k\geqslant0$,
\[
\mathcal{N}(k,k+3)\leqslant\mathcal{N}(a,a+3),
\qquad a=\mathcal{N}(1,3).
\]
Consequently,
\begin{equation}\label{eq:codim-three-stable}
\mathcal{N}(k,k+3)=\mathcal{N}(a,a+3)
\qquad(k\geqslant a).
\end{equation}
Using the computer-assisted equality $\mathcal{N}(1,3)=7$, this becomes
\[
\boxed{\mathcal{N}(k,k+3)=\mathcal{N}(7,10)\qquad(k\geqslant7).}
\]
\end{corollary}

\begin{proof}
Let $W=\sum_iU_i^\perp$ and $d=\dim W$. Theorem
\ref{thm:codim-three-normal} gives $d\leqslant a+3$. Since every
direction contains $W^\perp$, quotienting by this common direction
produces a configuration of affine $(d-3)$-planes in $\R^d$.
Adding a common factor embeds it into $\R^{a+3}$ as affine
$a$-planes. Thus
$\mathcal{N}(k,k+3)\leqslant\mathcal{N}(a,a+3)$.
For $k\geqslant a$, the reverse inequality follows by adding
a common $(k-a)$-dimensional factor. The numerical specialization
uses \cite[Theorem~1.1]{MRS} together with \cite{BLR}.
\end{proof}

\begin{remark}[Sharpness of the normal-span dimension]
The bound $a+3$ in Theorem~\ref{thm:codim-three-normal} is sharp.
Take a unit-distance family of $a$ lines $\ell_1,\ldots,\ell_a$
in $\R^3$, and work in $\R^3\oplus\R^a$. Let
\[
H_i=\{z\in\R^a|z_i=0\},\qquad L_i=\ell_i\times H_i.
\]
Each $L_i$ has codimension three, and
$d(L_i,L_j)^2=d(\ell_i,\ell_j)^2+d(H_i,H_j)^2=1$.
If $u_i$ is the direction of $\ell_i$, then
\[
U_i^\perp=u_i^\perp\oplus\R e_i.
\]
Since a parallel family of lines in $\R^3$ has at most three members,
a family of $a\geqslant7$ pairwise unit-distance lines cannot have all
directions parallel. Hence $\sum_i u_i^\perp=\R^3$.
Thus $\dim\sum_iU_i^\perp=a+3$.
\end{remark}

\subsection{A supermultiplicative sequence of stabilized values}

The fixed-codimension maxima $b_c$ themselves satisfy a product
inequality.

\begin{proposition}\label{prop:bc-supermultiplicative}
For all positive integers $c,d$,
\begin{equation}\label{eq:bc-supermultiplicative}
b_{c+d}\geqslant b_cb_d.
\end{equation}
In particular, $b_{c+1}\geqslant2b_c$ and
$b_{c+2}\geqslant a\,b_c$, where $a=\mathcal{N}(1,3)$.
\end{proposition}

\begin{proof}
By Theorem~\ref{thm:fixed-codim-stability}, each $b_c$ and $b_d$
is attained by a finite configuration. Choose configurations
realizing them, with parameters $(k_c,k_c+c)$ and $(k_d,k_d+d)$.
By adding common directions to either configuration, we may assume
$k_c\geqslant c$ and $k_d\geqslant d$. Thus both factors satisfy
the transversality condition in Theorem~\ref{thm:product}.
The transverse product gives a configuration of codimension $c+d$
with $b_cb_d$ members. Hence \eqref{eq:bc-supermultiplicative}.
Finally $b_1=2$ and $b_2=a$, by the endpoint formula and
codimension-two reduction.
\end{proof}

\section{Scope of the method and further questions}\label{sec:discussion}

\subsection{What is extended from H\"ofer's argument}
H\"ofer represents a three-dimensional line by $v=(l,m)\in\R^6$, with
$\norm l=1$ and $l\cdot m=0$, and lifts it to $vv\T$. For
\[
 M_1=\begin{pmatrix}I_3&0\\0&0\end{pmatrix},\qquad
 M_2=\begin{pmatrix}0&I_3\\I_3&0\end{pmatrix},
\]
his form is
\begin{equation}\label{eq:Hoefer-form}
 \B(A,B)=\tr(M_1AM_1B)+\tr(M_2AM_2B)
             -\tr(M_1A)\tr(M_1B).
\end{equation}
Lines give isotropic vectors, and nonparallel unit-distance pairs give
orthogonal vectors. A dimension bound $2r-3$ for a Pl\"ucker span of
dimension $r$, together with a minimal-dependence argument for rank-one
matrices, proves independence of the lifted representatives and then
the bound nine \cite[Definitions~1--2, Lemmas~4--5]{Hoefer}.

Here the higher-degree lift converts two polynomial kernels into common
total isotropy, but no analogous independence theorem is used.
A lower-rank locus can contain a continuous family: fixing a direction
projection $P$ and allowing $a\in\ker P$ to vary gives
$F_\rho=R_\rho=0$ for all pairs whenever $\rho>k$, irrespective of their
actual distance. The algebraic inverse image of the lifted span is
therefore partitioned into connected analytic pieces, and its rank is
forced downward on each piece. The actual unit-distance condition is
used again at the terminal parallel rank. This explains both why the
lifting is useful and why it is not by itself a finiteness proof.

The three-dimensional direction-independence lemma also fails in larger
ambient dimension. For example, the three lines
\[
 \R e_1,\qquad e_3+\R e_2,\qquad
 \tfrac12e_3+\tfrac{\sqrt3}{2}e_4+\R(e_1+e_2)
 \subset\R^4
\]
are pairwise at distance one although all directions lie in
$\Span(e_1,e_2)$. Their offsets form an equilateral triangle in the
orthogonal two-dimensional factor.

\subsection{Quantitative and structural questions}
\label{subsec:quantitative-questions}

The proof separates the cost of decomposing a bounded-degree real
algebraic set from the number of possible rank drops. Improving the
former without losing analyticity would improve the general upper
bound. The present double-exponential estimate is not expected to
reflect the scale of extremal configurations. For fixed $k$, the
amplification inequalities give linear lower growth in $n$ with
explicit coefficients, whereas for the balanced sequence
\[
S(c):=\mathcal N(c-1,2c-1),
\]
the transverse product construction gives exponential lower growth
in $c$. These are two different parameter regimes.

A particularly interesting problem concerns the optimal stabilization
threshold in fixed codimension. Recall that
\[
b_c:=\sup_{k\geqslant0}\mathcal N(k,k+c),
\]
which is finite by Theorem~\ref{thm:fixed-codim-stability}.
Define the minimal stabilization threshold by
\[
s_c:=\min\left\{s\geqslant0\big|
\mathcal N(k,k+c)=b_c
\text{ for every }k\geqslant s\right\}.
\]
The fixed-codimension stability theorem gives
\[
s_c\leqslant D_c-c,
\]
where $D_c$ is defined recursively by
\eqref{eq:Dc-recursion}. However, this bound is not expected to be
optimal. Indeed, in codimension two, the exact reduction
\[
\mathcal N(k,k+2)=\mathcal N(1,3)
\qquad(k\geqslant1)
\]
shows that $s_2=1$, whereas the general recursion gives
$D_2-2=2$.

Motivated by the codimension-one and codimension-two extremal
configurations, as well as the transverse product inequality, we
propose the following conjecture.

\begin{conjecture}[Optimal stabilization and extremal values]
\label{conj:optimal-stabilization}
For every integer $c\geqslant1$, the minimal stabilization threshold
and the stabilized extremal value are given by
\[
\boxed{
s_c=\left\lfloor\frac c2\right\rfloor,
\qquad
b_c=2^{\,c\bmod 2}\,
7^{\,\lfloor c/2\rfloor}.
}
\]
Equivalently,
\[
\boxed{
\mathcal N(k,k+c)
=
2^{\,c\bmod 2}\,
7^{\,\lfloor c/2\rfloor}
\qquad
\text{whenever }
k\geqslant\left\lfloor\frac c2\right\rfloor.
}
\]
\end{conjecture}

The numerical factor $7$ is motivated by the equality
$\mathcal N(1,3)=7$ established in \cite{MRS}, together with the
seven-cylinder construction of \cite{BLR}. The transverse product
theorem already yields the lower bound
\[
b_c\geqslant
2^{\,c\bmod 2}\,7^{\,\lfloor c/2\rfloor}.
\]
Thus the main new difficulty is to establish the matching upper
bound and the proposed optimal stabilization threshold.

A possible geometric explanation for
Conjecture~\ref{conj:optimal-stabilization} is an extremal
configuration decomposition principle: after passing to a suitable
realization and splitting off common directions, extremal
configurations of codimension $c$ should admit a decomposition
compatible with transverse products of
$\lfloor c/2\rfloor$ codimension-two configurations and, when $c$
is odd, one codimension-one configuration. Such a principle would
explain the multiplicative formula for $b_c$. It is important,
however, that the transverse product theorem only constructs
configurations with the predicted cardinality; it does not establish
that arbitrary extremal configurations possess such a decomposition.
Proving an appropriate rigidity or decomposition theorem is therefore
a central open problem.

Several further questions arise naturally within the framework
developed here:
\begin{itemize}
    \item Can the general stabilization threshold $D_c-c$ be
    substantially improved? More specifically, can one prove
    $s_c\leqslant\lfloor c/2\rfloor$ for every $c$, or establish
    the equality conjectured above?

    \item Can the stable extremal values $b_c$ be computed for
    $c\geqslant3$? In particular, is
    $b_3=14$, and more generally does
    \[
    b_c=2^{\,c\bmod 2}7^{\,\lfloor c/2\rfloor}
    \]
    hold for every positive integer $c$?

    \item Can one establish a geometric decomposition theorem for
    extremal configurations, showing that the transverse product
    construction captures all extremal configurations, possibly
    after an appropriate reduction of the ambient space?

    \item Can higher-rank restrictions, beyond
    $\dim(U_i+U_j)\leqslant k+1$, be reduced to bounded ambient
    dimension?

    \item Can one improve the transverse product constructions or
    find matching upper bounds for the balanced sequence $S(c)$?
\end{itemize}

Theorem~\ref{thm:top} provides a sharp model when all directions lie
in a common $(k+1)$-dimensional space, while
Example~\ref{ex:codim-three} shows why a general common-factor
assertion would be too strong. The conjectured extremal decomposition
would therefore require a more subtle structure than a universal
common-direction factorization.

All proofs of finiteness, parameter reduction, and product
constructions in this paper are independent of finite enumeration.
The numerical lower bounds using seven lines invoke the construction
in \cite{BLR}, and the exact specialization in
Corollary~\ref{cor:codim-seven} additionally invokes \cite{MRS}.
No global computer-free proof of $\mathcal N(1,3)=7$ is supplied.

\appendix
\section{Normalized affine Pl\"ucker data and quadratic recovery}\label{app:plucker}

Assume $1\leqslant k<n$. Choose an oriented orthonormal basis
$u_1,\ldots,u_k$ of a direction space $U$, and put
\[
 u=u_1\wedge\cdots\wedge u_k,\qquad \norm u=1,\qquad
 \eta=a\wedge u,
\]
where $a\in U^\perp$ is the canonical point of the affine plane.
The pair $(u,\eta)$ changes by a common sign under orientation reversal.
It is independent of the chosen base point along the plane, since
$(a+v)\wedge u=a\wedge u$ for $v\in U$.
Use the Euclidean inner products on exterior powers and let $\iota_v$
denote contraction by $v$.

\begin{proposition}\label{prop:Plucker-recovery}
The direction projection and canonical point satisfy
\begin{equation}\label{eq:Plucker-recovery}
 P_{ij}=\ip{\iota_{e_i}u}{\iota_{e_j}u},\qquad
 a_i=\ip{\iota_{e_i}\eta}{u}.
\end{equation}
Thus $(P,a)$ is a quadratic, orientation-independent expression in the
normalized affine Pl\"ucker data $(u,\eta)$.
\end{proposition}
\begin{proof}
The signed $(k-1)$-vectors obtained by omitting one factor from $u$ form an
orthonormal family. The contraction formula therefore gives
\[
 \ip{\iota_{e_i}u}{\iota_{e_j}u}
 =\sum_{\alpha=1}^k(e_i\cdot u_\alpha)(e_j\cdot u_\alpha)=P_{ij}.
\]
For the second identity, use
$\iota_{e_i}(a\wedge u)=a_i u-a\wedge\iota_{e_i}u$.
Since $a\perp U$, $\iota_a u=0$, and adjointness of exterior
multiplication and contraction gives
$\ip{a\wedge\iota_{e_i}u}{u}=\ip{\iota_{e_i}u}{\iota_a u}=0$.
Taking the inner product with the unit vector $u$ yields $a_i$.
Both expressions are unchanged by $(u,\eta)\mapsto(-u,-\eta)$.
\end{proof}

For $k=1,n=3$, the Hodge dual of $\eta=a\wedge u$ is the moment vector
$a\times u$ in H\"ofer's representation. Composing
\eqref{eq:Plucker-recovery} with \eqref{eq:veronese} gives a polynomial
lifting of degree at most $2(n+1)$ in these normalized data. This concerns
normalized coordinates, not a claim that normalization is a homogeneous
polynomial operation on projective space.

\end{document}